\documentclass[]{interact}

\usepackage{epstopdf}
\usepackage{subfigure}

\usepackage{natbib}
\bibpunct[, ]{(}{)}{;}{a}{}{,}

\theoremstyle{plain}
\newtheorem{theorem}{Theorem}[section]

\newtheorem{corollary}[theorem]{Corollary}
\newtheorem{proposition}[theorem]{Proposition}

\theoremstyle{definition}
\newtheorem{definition}[theorem]{Definition}

\theoremstyle{remark}

\theoremstyle{observation}
\newtheorem{observation}{Observation}
\newtheorem{conjecture}{Conjecture}

\usepackage{amsmath,amssymb}
\usepackage{url}
\usepackage[mathscr]{euscript}
\usepackage{stmaryrd}
\usepackage{mathabx}
\usepackage{dsfont}
\usepackage{graphicx}
\usepackage{xcolor,colortbl}
\newcommand {\e}{{\bf e}}
\newcommand{\scrA}{\mathscr{A}}
\newcommand{\scrB}{\mathscr{B}}

\newcommand{\scrL}{\mathscr{L}}
\newcommand{\scrM}{\mathscr{M}}
\newcommand{\scrP}{\mathscr{P}}
\newcommand{\bbR}{\mathbb{R}}

\newcommand{\En}{\operatorname{Em}}

\newcommand\defeq{\mathrel{\overset{\makebox[0pt]{\mbox{\normalfont\tiny\sffamily def}}}{=}}}
\newcommand{\Real}{\mathbb R}
\newcommand{\set}[1]{\left\{#1\right\}}                     
\newcommand{\abs}[1]{\left|#1\right|}                       
\newcommand{\bra}[1]{\left(#1\right)}
\newcommand{\brac}[1]{\left[#1\right]}

\usepackage{afterpage}

\begin{document}


\title{Strong mixed-integer formulations for the floor layout problem}

\author{
\name{Joey Huchette\textsuperscript{a}\thanks{CONTACT J. Huchette. Email: huchette@mit.edu}, Santanu S. Dey\textsuperscript{b}, and Juan Pablo Vielma\textsuperscript{c}}
\affil{\textsuperscript{a}Operations Research Center, Massachusetts Institute of Technology, Cambridge, MA\\ \textsuperscript{b}School of Industrial and Systems Engineering, Georgia Tech, Atlanta, GA\\ \textsuperscript{c}Sloan School of Management, Massachusetts Institute of Technology, Cambridge, MA}
}

\maketitle

\begin{abstract}
	The floor layout problem (FLP) tasks a designer with positioning a collection of rectangular boxes on a fixed floor in such a way that minimizes total communication costs between the components. While several mixed integer programming (MIP) formulations for this problem have been developed, it remains extremely challenging from a computational perspective.
	This work takes a systematic approach to constructing MIP formulations and valid inequalities for the FLP that unifies and recovers all known formulations for it. In addition, the approach yields new formulations that can provide a significant computational advantage and can solve previously unsolved instances. While the construction approach focuses on the FLP, it also exemplifies generic formulation techniques that should prove useful for broader classes of problems.
\end{abstract}



\maketitle

\section{Introduction}

The floor layout problem (FLP), also known as the (unequal areas) facility layout problem, is central to the design of objects such as factory floors and very-large-scale integration (VLSI) computer-chips. The designer is given a fixed rectangular floor and $N$ rectangular boxes to place onto the floor. Each box must sit completely on the floor, and they cannot overlap. Each box has a fixed area, but the widths and heights can be varied to change the shape, subject to constraints on the area and aspect ratio of the components. The objective is to minimize the weighted sum of the Manhattan norm distances between each pair of boxes.

The FLP can be naturally described as a disjunctive programming problem, which are often reformulated as mixed-integer programming (MIP) problem such as to take advantage of state-of-the-art MIP solvers. However, the FLP and its various MIP formulations have proven extremely difficult to solve to optimality. In this work, we take a systematic approach to generating MIP formulations for the FLP that unifies existing MIP formulations from the literature and leads to new formulations and valid inequalities. We also computationally compare the range of formulations, and show that the new approaches can be used to solve previously unsolved instances. The main contributions of this work include:

\begin{enumerate}
\item {\bf Case study on systematic construction of effective MIP formulations:} The number, heterogeneity, and complexity of existing MIP formulations for the FLP and the fact that it remains computationally challenging make it an excellent candidate for such a study. Through the use of the \emph{embedding formulation} approaches of \cite{Vielma:2015a} and through a systematic treatment of alternative disjunctive descriptions of the FLP, we are able to recover and unify all existing, seemingly ad-hoc, MIP formulations. In addition, we are able to derive new formulations that can provide a significant computational advantage and solve previously unsolved instances. While the study concentrates on specific characteristics of the FLP, it exemplifies generic formulation techniques and practices that should be useful for a wide range of problems.

\item {\bf Valid inequalities for alternative formulations of the FLP:} Using the embedding approach, we are able to construct a variety of new valid inequalities for FLP. One key of the embedding formulation approach of this work is the flexible use of $0/1$ variables to model disjunctive constraints or unions of polyhedra. However, such flexibility can cloud the ``interprebility'' of the $0/1$ variables, which is often needed to construct valid inequalities to strengthen formulations. In this work, we show how to construct and translate valid inequalities between formulations of the FLP in such a way that allows us to state a broad class of valid inequalities in a generic form.

\item {\bf Comprehensive computational study of FLP:} While the FLP has been extensively studied, most existing works compare only a small subset of formulations and valid inequalities when making comparison. In this work we collect several instances from the literature to construct a publicly available library and present a comprehensive computational study of existing and new formulations on this library. Furthermore, our systematic approach allows us to compare a host of formulations and a wide range of common valid inequalities when making our comparison. In particular, while no single formulation seems to be dominant, we may offer a small collection of techniques which prove particularly effective for the FLP. Furthermore, we also study various theoretical and practical aspects of these approaches that help explain their success.

\end{enumerate}

The remainder of this work is organized as follows. In Section~\ref{sec:lit} we present a literature review of the existing solution techniques for the FLP. In Section~\ref{sec:def} we formally define the FLP and show how it can be cast as a disjunctive programming problem. In Section~\ref{sec:formulations} we review the formulation techniques we use to transform the FLP into a MIP, and in Section~\ref{sec:twobox} we use the techniques to construct formulations that are based on the interaction of two boxes at a time. Then in Section~\ref{sec:inequalities} we develop valid inequalities that can be used to strengthen formulations and show how they can be translated from one formulation to the other. In Section~\ref{sec:inequalities} we also restrict attention to the interaction of two boxes at a time, so in Section~\ref{sec:multibox} we develop formulations and inequalities that are based on the interaction of larger collections of boxes. Finally, in Section~\ref{sec:computations} we present results of our computational experiments, and in Section~\ref{sec:conclusions} we present a brief summary of this work. Complementary material and omitted proofs are included in the Appendix.

\section{Literature review}\label{sec:lit}
The floor layout problem can be viewed as a specific version of a general layout problem that consists of orthogonally packing rectangular pieces onto a rectangular floor; \cite{Drira:2007} offer a taxonomy of variations of the FLP and its relatives. Originally studied primarily in the context of factory design, the emergence of the field of very-large scale integration (VLSI) computer-chip design saw renewed interest in layout problems such as the FLP.


Broadly, algorithmic approaches to these layout problems can be grouped into two classes: exact and heuristic. Exact algorithms were predominant in the earlier literature, although the boom of applications in computer-chip design require solving large scale instances beyond the reach of existing exact approaches. As a result, a bevy of work has appeared over the past three decades, proposing heuristic approaches to produce good solutions for large-scale instances. Much of the work applies existing metaheuristic frameworks to the FLP, for example \cite{Meller:1996a} and \cite{Tate:1995}. Contrastingly, many of the novel heuristics for the FLP take advantage of ideas and machinery from mathematical programming: e.g. \cite{Camp:1991,Anjos:2006,Bernardi:2010,Bernardi:2013,Jankovits:2011,Lin:2011,Luo:2008,Luo:2008a}, and \cite{Liu:2007}, albeit in a way that cannot prove optimality. We note in particular the surveys of \cite{Meller:1996} and \cite{Singh:2006}, which collect pointers to much of the heuristic literature.

In keeping with the MIP approach taken in this paper, we will survey the existing exact methods for the FLP in detail. Early work can be traced back to \cite{Bazaraa:1975}, who studies a discretized version of the FLP. \cite{Meller:1999} introduced a natural MIP model for the FLP, along with a collection of valid inequalities and techniques to help reduce solution time. \cite{Sherali:2003} introduces novel formulations for a single pair of boxes, as well as useful computational techniques such as symmetry breaking constraints and branching priorities. \cite{Castillo:2005} presents a new MIP formulation for the FLP with fewer binary variables, alongside a number of additional formulations and approaches inspired by nonlinear and mixed-integer nonlinear optimization. \cite{Meller:2007} presents another formulation inspired by a technique from \cite{Murata:1996} that reduces redundancy in the solution set. As detailed in the following section, the inclusion of certain non-linear area constraints in the FLP result on its formulations being second-order-cone MIP (MISOCP) problems. Given that early formulations were developed before the availability of efficient MISOCP solvers, careful attention has been paid on constructing and proving desirable properties for specific linear approximations for the nonlinear area constraints in \cite{Sherali:2003} and \cite{Castillo:2005a}.

The FLP has a natural one-dimensional analogue in the single-row floor (facility) layout problem, which asks for an optimal layout of $N$ boxes of fixed length in a straight line. This problem is already NP-hard, and strong formulations and cutting planes have been developed for the problem by \cite{Amaral:2006,Amaral:2008} and \cite{Amaral:2012}.


An intriguing line of research has investigated the FLP from the dual perspective, attempting to construct tight lower bounds. This is of particular interest for the FLP, where relaxations typically give poor bounds, even with strengthening valid inequalities. \cite{Amaral:2009} presents a lower bounding technique for the single-row FLP. Another line of work investigates using semidefinite programming formulations to construct bounds for the FLP by \cite{Takouda:2005} and the single-row FLP by \cite{Anjos:2005}. \cite{Anjos:2008} leverage the semidefinite approach to produce optimal solutions for the single-row FLP using a cutting-plane approach, and to produce high-quality solutions for larger instances in \cite{Anjos:2009}.  \cite{Huchette:2015a} present a combinatorial dual bounding scheme for the FLP and compare it against existing techniques.

\section{Preliminaries: Notation, defining constraints, and the disjunctive formulation}\label{sec:def}
Consider a rectangular floor $[0,L^x] \times [0,L^y]$ for $L^x,\,L^y>0$. There is a collection of $N$ boxes $\{\scrB_i\}_{i=1}^N$ to place on the floor, each with a target area $\alpha_i>0$ and maximum allowed aspect ratio $\beta_i>0$. Denote the set of all pairs with $\scrP = \{(i,j) \in \llbracket N \rrbracket^2 : i < j\}$ where $\llbracket N \rrbracket \defeq \{1,\ldots,N\}$. With each pair of boxes $(i,j) \in \scrP$, there is an associated nonnegative unit communication cost $p_{i,j}$. The floor layout problem then is to optimally lay out each box completely on the floor, such that the area and aspect ratio constraints are satisfied, and such that no two boxes overlap.

Natural decision variables for each box $\scrB_i$ are the position of its center $(c^x_i,c^y_i)$ and the lengths in each direction $(\ell^x_i,\ell^y_i)$. The objective function used is based on the so-called ``Manhattan'' norm:
\begin{equation} \label{eqn:objective}
	\sum_{(i,j) \in \scrP}p_{i,j}\left( |c^x_i-c^x_j| + |c^y_i-c^y_j| \right).
\end{equation}

Most of the constraints described are simple to describe with linear or conic inequalities. For instance, $\scrB_i$ lies completely on the floor iff
\begin{equation}
	\frac{1}{2}\ell^s_i \leq c^s_i \leq L^s - \frac{1}{2}\ell^s_i \quad \forall s \in \{x,y\}, i \in \llbracket N \rrbracket. \label{eqn:sitb}
\end{equation}
The area constraints take the form
\begin{equation} \label{eqn:area}
	\ell^x_i \ell^y_i \geq \alpha_i \quad \forall i \in \llbracket N \rrbracket,
\end{equation}
which is second-order-cone-representable \citep{Alizadeh:2003}. 
The aspect ratio constraints take the form
\begin{equation}
	\max\left\{\frac{\ell^x_i}{\ell^y_i}, \frac{\ell^y_i}{\ell^x_i}\right\} \leq \beta_i \quad \forall i \in \llbracket N \rrbracket.
\end{equation}
This can be represented with two linear constraints per box, but it can also be enforced on the FLP merely through bounds on the widths of the boxes.
\begin{observation}[Section 2.3 of \cite{Castillo:2005}]
Along with the area constraints, imposing the following bounds on the box widths is sufficient to impose the aspect ratio constraints:
\begin{subequations} \label{eqn:bounds}
\begin{alignat}{2}
	\ell^s_i &\leq ub^s_i \defeq \min\left\{ \sqrt{\alpha_i\beta_i}, L^s \right\} \quad& \forall s \in \{x,y\}, i \in \llbracket N \rrbracket \label{eqn:upperbound} \\
	\ell^s_i &\geq lb^s_i \defeq \frac{\beta_i}{ub^s_i} \quad& \forall s \in \{x,y\}, i \in \llbracket N \rrbracket\label{eqn:lowerbound}
\end{alignat}
\end{subequations}
\end{observation}
We note that, since $\beta_i > 0$, we have that $lb^s_i > 0$ for each $s \in \{x,y\}$ and $i \in \llbracket N \rrbracket$.

 \subsection{Disjunctive formulation for the FLP}
 The last remaining constraint for the FLP requires that the boxes cannot overlap on the floor. One natural way to formulate this is by requiring each pair $\scrB_i$ and $\scrB_j$ to be separated in either the $x$ direction or the $y$ direction (or both).

\begin{definition} \label{defn:precede}
 We say that $\scrB_i$ \emph{precedes} $\scrB_j$ in direction $s$ (denoted by $\scrB_i \leftarrow_s \scrB_j$) if
\begin{equation} \label{eqn:precedes}
	c^s_i + \frac{1}{2}\ell^s_i \leq c^s_j - \frac{1}{2}\ell^s_j.
\end{equation}
\end{definition}
Therefore, we can enforce the constraint that $\scrB_i$ and $\scrB_j$ do not overlap with the disjunctive constraint $D^4_{i,j} \defeq \bigvee_{k=1}^4 d^k_{i,j}$,
where
\[
    d^1_{i,j} = \scrB_i \leftarrow_y \scrB_j, \quad\quad d^2_{i,j} = \scrB_i \leftarrow_x \scrB_j, \quad\quad d^3_{i,j} = \scrB_j \leftarrow_y \scrB_i, \quad\quad d^4_{i,j} = \scrB_j \leftarrow_x \scrB_i.
\]
We omit the subscript and use $D^4$ when the meaning is clear from context and we refer to each $d^k$ as a \emph{branch} of the disjunction.

\begin{figure}
\centering
\includegraphics[width=.32\linewidth]{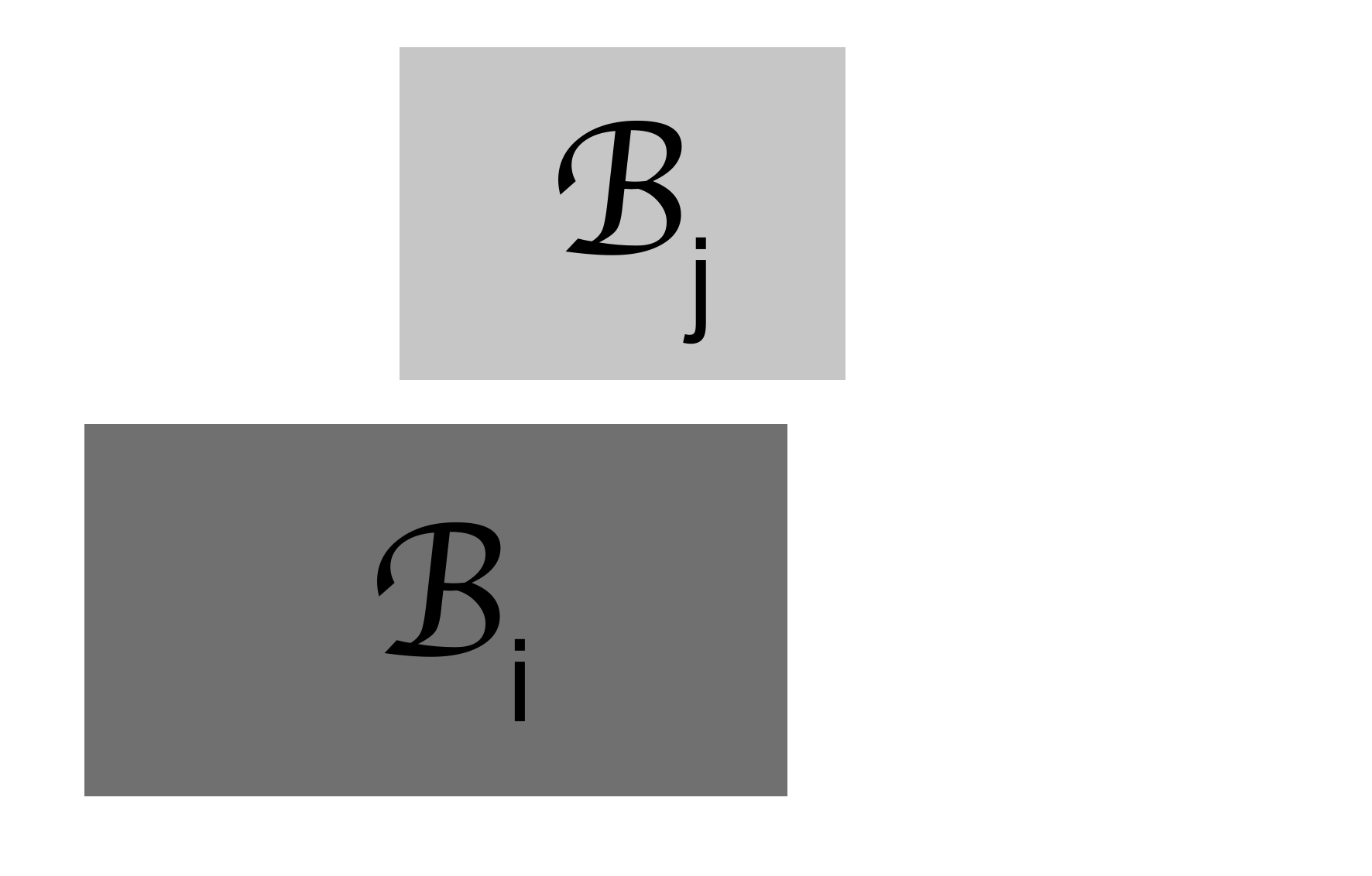}
\includegraphics[width=.32\linewidth]{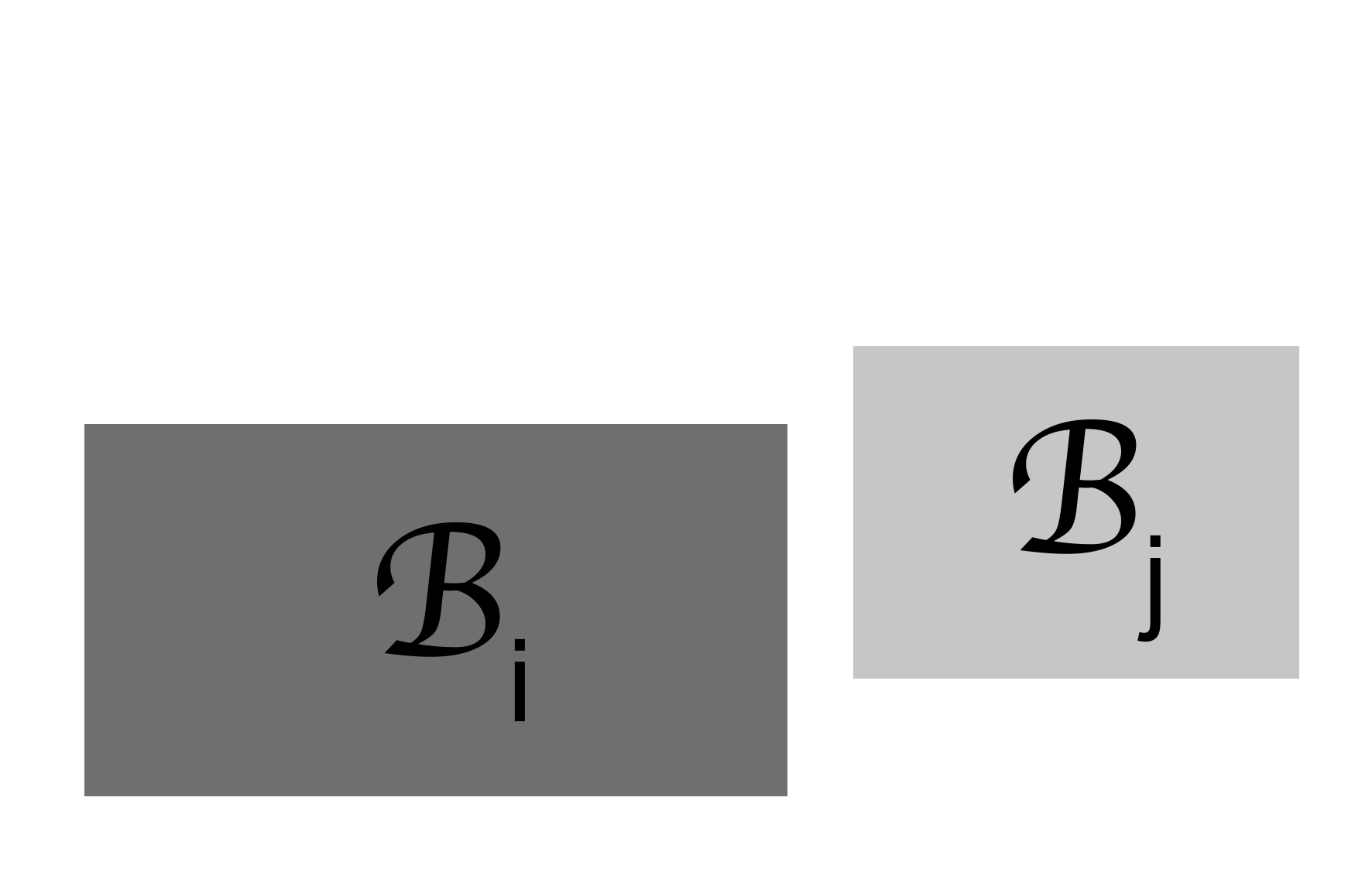}
\includegraphics[width=.32\linewidth]{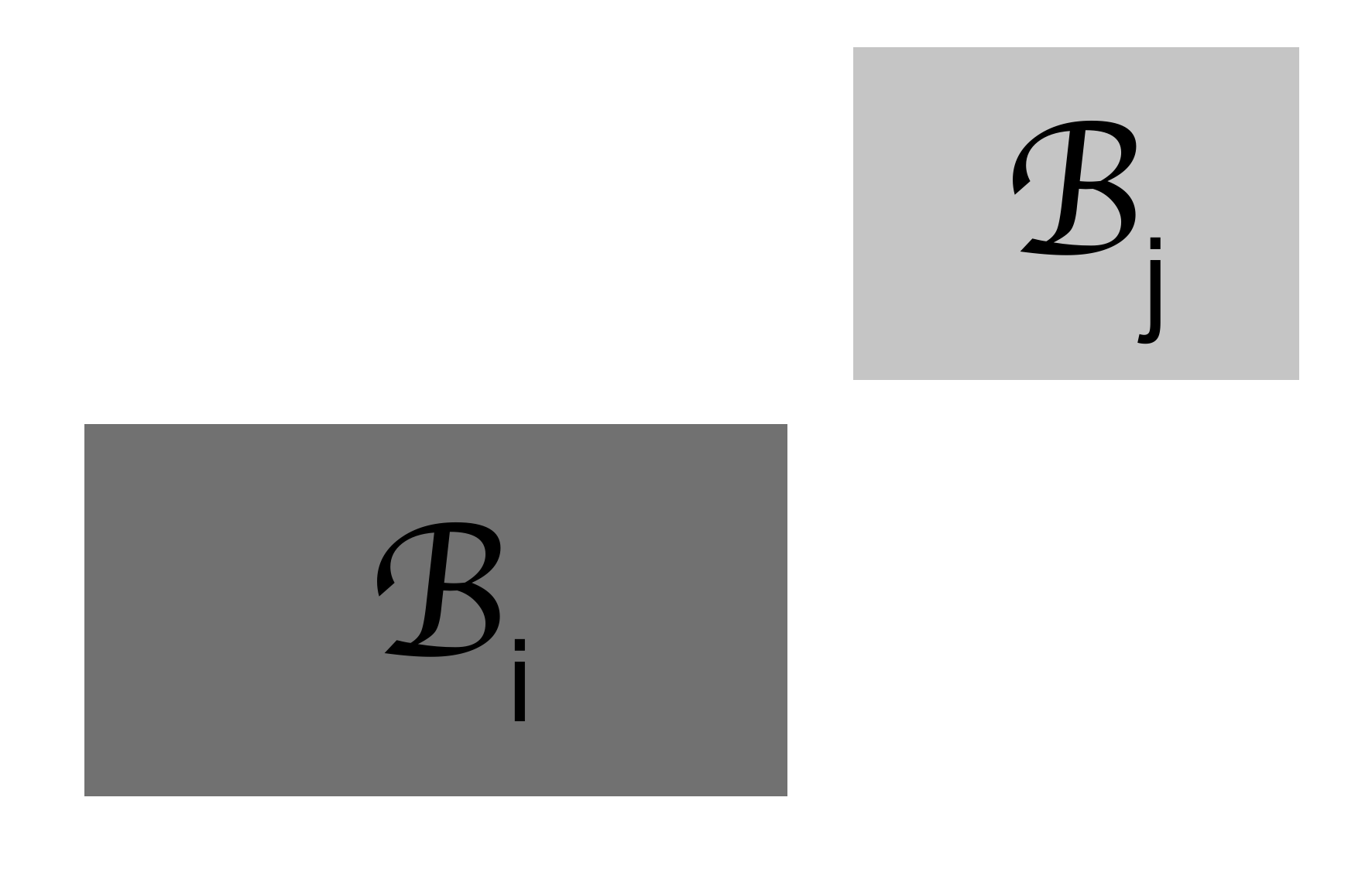}
\caption{Two boxes whose layout satisfies $d^1_{i,j}$ (Left), $d^2_{i,j}$ (Center), and both $d^1_{i,j}$ and $d^2_{i,j}$ (Right).}
\label{fig:two-boxes}
\end{figure}

Then the set of all feasible layouts is given by the disjunctive set
\begin{equation} \label{eqn:set-of-feasible-layouts}
	\scrL \defeq \left\{ (c,\ell) \in \bbR^{2N+2N}: \eqref{eqn:sitb},\eqref{eqn:bounds},\eqref{eqn:area}, \bigwedge_{(i,j) \in \scrP} D^4_{i,j} \right\},
\end{equation}
and a (nonlinear) disjunctive programming formulation of the FLP is given by $ \min\set{\eqref{eqn:objective}\,:\,   (c,\ell) \in \scrL }$. The main objective of this work is to transform this disjunctive programming problem into a mixed-integer formulation that can be solved by  off-the-shelf optimization software. To achieve this, we first focus our attention to a single pair of boxes $\scrB_i$ and $\scrB_j$ for $(i,j) \in \scrP$, in which case we are interested in constructing a MIP formulation of $\scrL_{i,j} \defeq \left\{ (c_i,c_j,\ell_i,\ell_j) \in \bbR^8: \eqref{eqn:sitb},\;\eqref{eqn:bounds},\;\eqref{eqn:area},\; D^4 \right\}$
If we omit the nonlinear area constraints \eqref{eqn:area} we obtain the set  $\hat{\scrL}_{i,j} \defeq \left\{ (c_i,c_j,\ell_i,\ell_j) \in \bbR^8: \eqref{eqn:sitb},\;\eqref{eqn:bounds},\; D^4 \right\}$,
which can be written as the union of (bounded) polyhedra $\hat{\scrL}_{i,j} = \bigcup_{k=1}^4 P^k_{i,j}$, where $P^k_{i,j} = \left\{ (c_i,c_j,\ell_i,\ell_j) \in \bbR^8 : \eqref{eqn:sitb},\;\eqref{eqn:bounds},\; d^k_{i,j} \right\}$. We can then use various techniques to construct a linear MIP formulation of $\hat{\scrL}_{i,j}$, which can be combined with the area constraints to obtain a second-order-cone MIP (MISOCP)
formulation of $\scrL_{i,j}$. Finally, we can combine these formulations for all pairs $(i,j) \in \scrP$ and an appropriate linearization of the objective function \eqref{eqn:objective} to obtain a MISOCP of the complete problem. Much of our analysis will focus on constructing such MIP formulations of $\hat{\scrL}_{i,j}$, which we denote the \emph{pairwise FLP}. However, in Section \ref{sec:objective} and Section \ref{ss:multibox-cuts} we consider the strengthening of the final formulation by explicitly considering the objective function and larger collections of boxes, respectively, when formulating the disjunctive constraint.


\section{Constructing MIP formulations for disjunctive constraints}\label{sec:formulations}

MIP formulations for unions of  polyhedra such as $\hat{\scrL}_{i,j}$ can be divided into \emph{extended formulations} that uses both continuous and $0/1$ auxiliary variables \citep[Section 5]{Vielma:2015} and \emph{non-extended (projected) formulations} that only use the  $0/1$ auxiliary variables that are strictly necessary to build a valid formulation \citep[Section 6]{Vielma:2015}. Note that, in general, it is not possible to construct a MIP formulation for unions of polyhedra in the original space, so some ($0/1$) variables are needed for representability. Standard extended formulations by \cite{Balas:1998} and \cite{Jeroslow:1984} have the desirable property that their Linear Programming (LP) relaxations have extreme points that naturally satisfy the integrality requirements on the $0/1$ auxiliary variables; we  call such formulations \emph{integral} or \emph{ideal}. In contrast, non-extended formulations often fail to be ideal, but can be much smaller. For instance, the following proposition shows a non-extended formulation of $\hat{\scrL}_{i,j}$ obtained through the classical \emph{big-$M$} approach.
\begin{proposition} \label{observation:bigM}
	The following is a formulation for $\hat{\scrL}_{i,j}$:
	\begin{subequations}\label{eqn:unary-bigM}
	\begin{align}
		\frac{1}{2}\ell^s_k \leq c^s_k \leq L^s - \frac{1}{2}\ell^s_k \quad &\forall s \in \{x,y\}, k \in \{i,j\} \\
		c^y_i + \frac{1}{2}\ell^y_i \leq c^y_j - \frac{1}{2}\ell^y_j + L^y(1-v_1), \quad\quad&
		c^x_i + \frac{1}{2}\ell^x_i \leq c^x_j - \frac{1}{2}\ell^x_j + L^x(1-v_2) \\
		c^y_j + \frac{1}{2}\ell^y_j \leq c^y_i - \frac{1}{2}\ell^y_i + L^y(1-v_3), \quad\quad&
		c^x_j + \frac{1}{2}\ell^x_j \leq c^x_i - \frac{1}{2}\ell^x_i + L^x(1-v_4) \\
		lb^s_k \leq \ell^s_k \leq ub^s_k \quad &\forall s \in \{x,y\}, k \in \{i,j\}\\
	    \sum_{i=1}^4 v_i = 1, \quad\quad& v \in \{0,1\}^4.
	\end{align}
	\end{subequations}
\end{proposition}
Formulation \eqref{eqn:unary-bigM} only uses four $0/1$ auxiliary variables (and no continuous auxiliary variables) and is about four times smaller than the standard ideal extended formulation for  $\hat{\scrL}_{i,j}$ (see Appendix~\ref{extendedformulationsec}). While formulation \eqref{eqn:unary-bigM} is not guaranteed to be ideal, its smaller size can still result in a computational advantage over the ideal extended formulation. In addition, using various techniques it is sometimes possible to strengthen non-extended formulations considerably. For this reason, we concentrate on constructing non-extended formulations for the FLP by using three techniques: (1) the flexible use of $0/1$ auxiliary variables provided by the  \emph{embedding formulations} approach, (2) alternative definitions of the disjunctive constraints $D^4_{i,j}$, and (3) the consideration of various common linear inequalities when building the disjunctions. We now provide some simple examples of applying these techniques to build MIP formulations for disjunctive sets.

\subsection{Selecting encodings and alternative MIP formulations: The embedding approach}

We begin our description of the embedding formulation approach of \cite{Vielma:2015a} by re-interpreting  \eqref{eqn:unary-bigM} as a formulation for the \emph{embedding} of $\hat{\scrL}_{i,j}$ in a higher dimensional space. Indeed,  $(c_i,c_j,\ell_i,\ell_j,v)$ is feasible for  \eqref{eqn:unary-bigM} if and only if it belongs to the embedding  of $\hat{\scrL}_{i,j} = \bigcup_{k=1}^4 P^k_{i,j}\subseteq \Real^8$  into $\Real^8\times \set{0,1}^4$  given by
\begin{equation} \label{eqn:bigM-set}
	\bigcup_{k=1}^4 P^k_{i,j} \times \{\e^k\},
\end{equation}
where $\e^k$ is the $k$-th unit vector ($\e^k_l=0$ for $l\neq k$ and $\e^k_k=1$). We say that representation \eqref{eqn:bigM-set} \emph{embeds} $\hat{\scrL}_{i,j}$, which lives in the space of the $(c_i,c_j,\ell_i,\ell_j)$ variables, into the space of the  $(c_i,c_j,\ell_i,\ell_j,v)$ variables. It achieves this by pairing each  of the four polyhedra  $ P^k_{i,j}$ with a unique binary vector $\e^k$, and so \emph{encoding} the disjunctive constraint. Any valid formulation for \eqref{eqn:bigM-set} implies a valid formulation for $\hat{\scrL}_{i,j}$, since $\hat{\scrL}_{i,j}$ is the orthogonal projection of  $ \bigcup_{k=1}^4 \left(P^k_{i,j} \times \{\e^k\}\right) $ onto the $(c_i,c_j,\ell_i,\ell_j)$ variables. However, representation \eqref{eqn:bigM-set} also makes  explicit the role of the $v$ variables: $v=\e^k$ implies that $(c_i,c_j,\ell_i,\ell_j)\in P^k_{i,j}$. In other words, the possible values $\set{\e^k}_{k=1}^4$ of $v$ \emph{encode} the selection among the polytopes $ P^k_{i,j}$.

The key for the flexibility of the embedding approach is noting that this encoding can use any family of pairwise distinct $0/1$ vectors in place of the unit vectors $\e^k$. The following definition formalizes this approach in our context, where we explicitly separate the disjunctive constraint $D$ from the common constraints $Q$, which must be satisfied by all branches of the disjunction.


\begin{definition}\label{embeddingdef}
    Take a polyhedra $Q \subseteq \bbR^d$, a disjunctive constraint $D = \bigvee_{k=1}^K[A^kx \leq b^k]$, and an encoding $C \defeq \{h^k\}_{k=1}^K \subseteq \{0,1\}^r$ of pairwise distinct vectors. A \emph{non-extended (linear) MIP formulation} for $\{x \in Q : D\}$ is any (linear) MIP formulation for the embedding 
    \[
        \En(Q,D,C) \defeq \bigcup\nolimits_{k=1}^K \left\{x \in Q : A^k x \leq b^k\right\} \times \{h^k\}
    \]
    that uses only $d+r$ variables.
\end{definition}

Standard formulation approaches are recovered when choosing the unit vectors $U^K \defeq \{\e^k\}_{k=1}^K \subseteq \{0,1\}^K$, which we denote the \emph{unary encoding}, as it uses one bit per branch of the disjunction. For example, we have that \eqref{eqn:unary-bigM} is a formulation for $\En(Q, D^4, U^4)$ with $Q = \left\{ (c_i,c_j,\ell_i,\ell_j) \in \bbR^8 : \eqref{eqn:sitb},\;\eqref{eqn:bounds} \right\}$. However, the real flexibility comes from the possibility of non-unary encodings, as the specific assignment of codes to branches of the disjunctions does not change the structure of the formulation. For instance, to obtain a valid formulation for $\En(Q, D^4, \tilde{U}^4)$ with $\tilde{U}^4=\set{\e^2,\e^1,\e^3,\e^4}$ we simply need to interchange $v_1$ and $v_2$ in \eqref{eqn:unary-bigM}. In contrast, for other types of encodings the specific assignment can be significant in terms of the complexity of the resulting embedding object and formulations (e.g. see  Section~\ref{sec:binaryform} and \cite{Vielma:2015a}).

Deriving ideal non-extended formulations for embeddings with any encoding can be done using a geometric construction introduced in
\cite{Vielma:2015a}. However, such construction can be hard to analyze, and many choices of encodings may naturally have very large ideal formulations (i.e. many inequalities). Fortunately, non-extended formulations can also be constructed using ad-hoc approaches or through simple constructions such as a generalization of the big-$M$ approach to arbitrary encodings. In the coming sections we will see how this generic approach can be used to construct a range of formulations for our disjunctive set. In particular, varying the ingredients $Q$, $D$, and $C$ lead to different embedding objects $\En(Q,D,C)$, which in turn will necessitate different formulations. In the following subsections, we provide such examples for varying inputs $D$ and $Q$.

\subsection{Alternative disjunctive formulations and common constraints}
\subsubsection{Alternative logical representations ($D$)}
The new formulation for the FLP proposed in Section \ref{sec:refined-disjunction} hinges on a logical refinement of the disjunction $D^4$ that removes many redundant solutions from the resulting formulation.

\begin{figure}[htpb]
  \begin{center}
  \subfigure[$\mathcal{M}_1$]{\label{fig1}\includegraphics[scale=0.47]{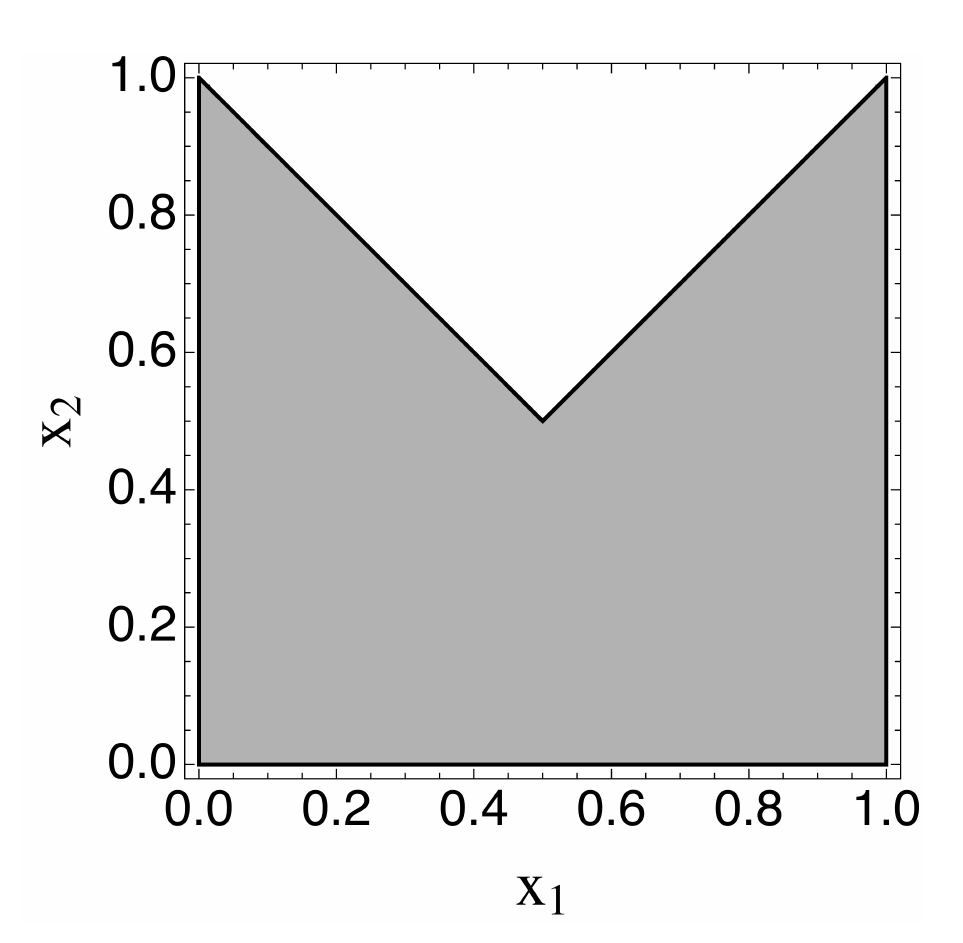}}
  \subfigure[$\mathcal{M}_2$ and $\{x \in Q_2^A: D_2\}$]{\label{fig3}\includegraphics[scale=0.47]{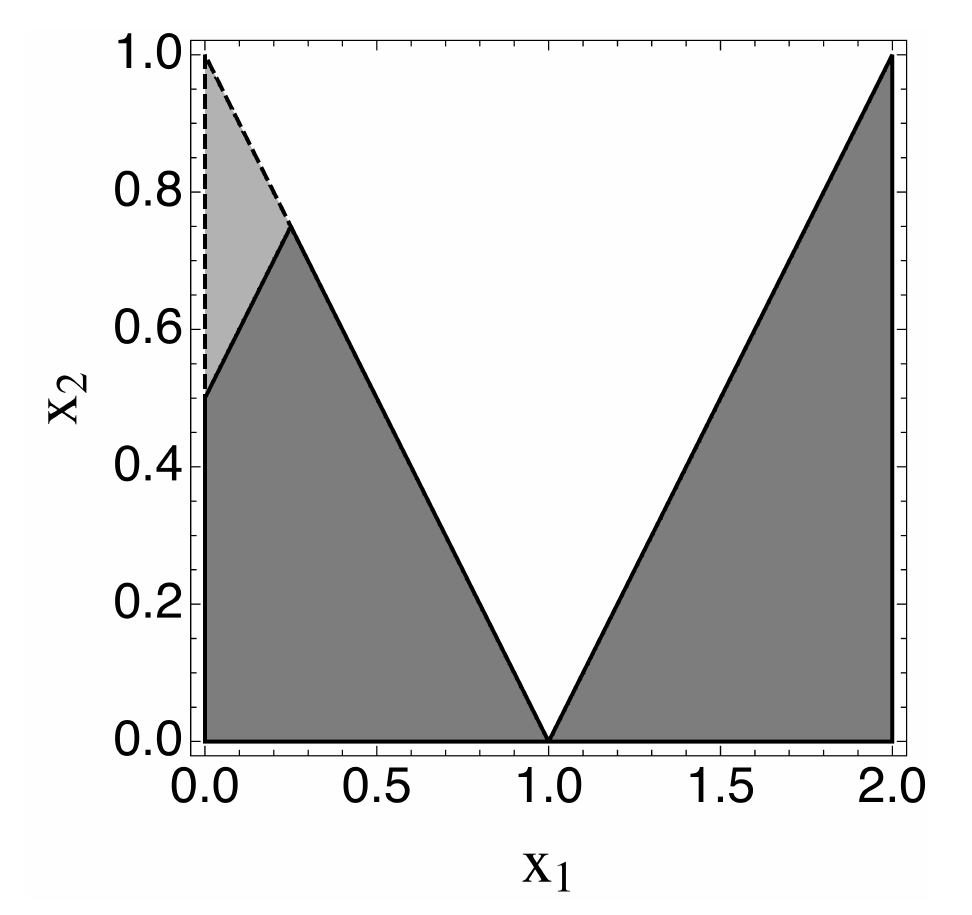}}
          \subfigure[$\mathcal{M}_3$]{\label{fig4}\includegraphics[scale=0.47]{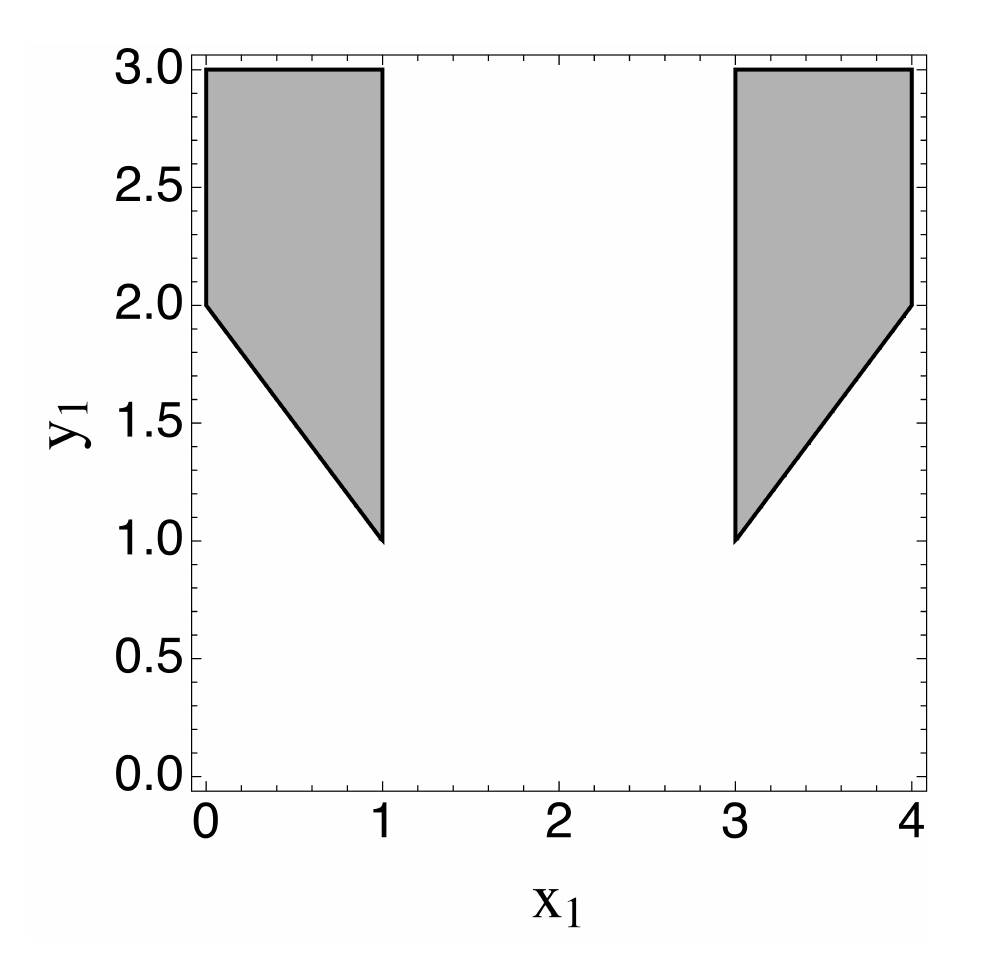}}
  \end{center}
  \caption{Disjunctive constraints.}\label{pwl2fig}
  \end{figure}

To illustrate this idea, we provide a simple example, which is independent of the FLP. Consider the disjunctive constraint
$D_1^A \defeq \brac{x_1+x_2\leq 1} \vee \brac{  x_2\leq x_1}$
and common linear constraints $Q_1 \defeq \set{x\in \Real^2\,:\, 0\leq x_i \leq 1 \quad \forall i\in \llbracket 2 \rrbracket }$, for which $\mathcal{M}_1\defeq\{x \in Q_1: D_1^A\}$ is depicted in Figure~\ref{fig1}. Because the two alternatives of $D^A_1$ intersect, we can define an alternative disjunction
\[
	D^B_1 \defeq \brac{ x_1+x_2\leq 1, x_1\leq x_2} \vee \brac{ x_1+x_2\leq 1, x_2\leq x_1}\vee \brac{x_1+x_2\geq 1, x_2\leq x_1},
\]
for which $\{x \in Q_1: D^A_1\}=\{x \in Q_1: D^B_1\}$. Using $D^B$ instead of $D^A$ could lead to larger formulations, since more branches on the disjunction will require longer codes to satisfy the distinctness property. However, it also reduces redundancy or symmetry, phenomena which are known to reduce the effectiveness of mixed-integer solvers. In particular, we note that the point $(1/2,1/4)$ satisfies both branches of $D^A_1$, but only one branch of $D^B_1$. Using the embedding approach for some encoding $C$, this gives to two feasible points in $\En(Q_1,D^A_1,C)$ which correspond to $(1/2,1/4)$ and differ only in their assigned code.
  
\subsubsection{Adding common constraints ($Q$)} \label{sec:common-constraints}
One advantage of the embedding framework as described in Definition~\ref{embeddingdef} is that it allows us to exploit the fact that the FLP has many constraints beyond the disjunctive constraint $D^4$. In particular, we can pick which constraints we include in the ground set $Q$ to combine with the disjunctive constraint $D$ and build a formulation, and which constraints are added after the formulation process. Such choice can significantly change the strength of the final formulation.

To illustrate this effect consider the following simple example, which is independent of the FLP. Take the disjunctive constraint  $D_2 \defeq \brac{ x_1+x_2\leq 1} \vee \brac{   1+x_2\leq x_1}$  and the set of linear inequalities
\begin{subequations}\label{alllinearineq}
\begin{alignat}{3}
0 &\leq x_i \leq 2\label{linearineq1}&\quad& i\in \llbracket 2 \rrbracket\\
x_2 &\leq x_1 + \frac{1}{2}.\label{linearineq}
\end{alignat}
\end{subequations}
Suppose we want to construct  a MIP formulation for $\mathcal{M}_2\defeq\set{x\in \Real^2\,:\,\eqref{alllinearineq}, D_2}$ depicted  by the dark shaded region in Figure~\ref{fig3}. If we let $Q_2^A  \defeq  \set{x\in \Real^2\,:\, \eqref{linearineq1} }$ we have that both $P^1=\set{x\in Q_2^A\,:\, x_2\leq 1-x_1}$ and $P^2=\set{x\in Q_2^A\,:\, 1+x_2\leq x_1}$ are bounded and hence satisfy the conditions of Definition~\ref{embeddingdef}. Then, we can at first ignore linear inequality \eqref{linearineq} and construct a formulation for $\{x \in Q_2^A: D_2\}=P^1\cup P^2$ (depicted by the and light shaded region in Figure~\ref{fig3}) and then impose \eqref{linearineq} on the resulting formulation. For instance, an ideal formulation of $ \En\bra{Q_2^A,D_2,\set{\e^1,\e^2}}$ is given by
 \begin{subequations}\label{partialform}
\begin{alignat}{3}
0\leq x_2&\leq 3-2v_1- x_1\\
1-v_1 \leq x_1&\leq 2-v_1\\
1+x_2&\leq x_1+2v_1\label{partialformreplace}\\
v_1+v_2&=1\\
v&\in \set{0,1}^2.
\end{alignat}
\end{subequations}
A formulation of $\mathcal{M}_2$ is then given by \eqref{partialform} and \eqref{linearineq}. However, a second option is to include all inequalities into $Q_2^B  \defeq  \set{x\in \Real^2\,:\, \eqref{linearineq1} \text{--}\eqref{linearineq}}$ and directly construct a formulation of $\{x \in Q_2^B: D_2\}=\mathcal{M}_2$. For instance, an ideal formulation of $ \En\bra{Q_2^B,D_2,\set{\e^1,\e^2}}$ is given by \eqref{partialform} with \eqref{partialformreplace} strengthened to
\begin{equation}\label{linearineqwithip}
1+x_2\leq x_1+\frac{3}{2}v_1.
\end{equation}
We can check that $\bra{x_1,x_2,v_1,v_2}=\bra{1,1,\frac{1}{2},\frac{1}{2}}$ is feasible for the LP relaxation of \eqref{partialform} and \eqref{linearineq}, but it does not satisfy \eqref{linearineqwithip}. Hence, the formulation obtained by considering all common linear inequalities is stronger that the one obtained by first ignoring \eqref{linearineq}.

A similar strengthening effect can occur when auxiliary variables and linear inequalities used to model other aspects of a mathematical programming problem are included in the common constraints. This is the case with the FLP and the nonlinear objective function \eqref{eqn:objective}, which may be linearized with auxiliary variables and constraints.

For a simple motivating example, let $ D_3 \defeq \brac{x_1\leq 1} \vee \brac{ x_1 \geq 3}$, $ Q_3^A \defeq \set{x_1\in \Real\,:\, 0\leq x_1\leq 4}$ and $\mathcal{M}_3^A\defeq\set{x_1\in Q_3^A\,:\, D_3}$, and suppose we want to solve $\min\set{\abs{x_1-2}\,:\, x_1\in \mathcal{M}_3^A}$. An ideal formulation for  $ \En\bra{Q_3^A,D_3,\set{\e^1,\e^2}}$ is given by
\begin{equation} \label{formulsimpleminus}
    3-3v_1\leq x_1 \leq 4-3v_1, \quad\quad v_1+v_2=1, \quad\quad v \in \{0,1\}^2,
 \end{equation}
 which together with a standard LP modeling trick to linearize the absolute value in the objective leads to the MIP formulation of the complete problem given by
 \begin{equation}\label{formulsimple}
 \min\set{y_1\,:\, x_1 -2 \leq y_1,\quad -x_1 +2 \leq y_1,\quad  \eqref{formulsimpleminus}}.
 \end{equation}
Alternatively, we could instead include the linearization trick in the common constraints to obtain $ Q_3^B \defeq \set{\bra{x_1,y_1}\in \Real^2\,:\, 0\leq x_1\leq 4,\quad x -2 \leq y_1,\quad -x_1 +2 \leq y_1 }$ and $\mathcal{M}_3^B\defeq\set{\bra{x_1,y_1}\in Q_3^B\,:\, D_3}$, depicted in Figure~\ref{fig4}. An integral formulation for $ \En\bra{Q_3^B,D_3,\set{\e^1,\e^2}}$ is given by \eqref{formulsimpleminus} plus
\begin{equation} \label{formulsimpletwo}
 x_1 -2 +2v_1 \leq y_1, \quad\quad -x_1 +2 +2(1-v_1) \leq y_1,
\end{equation}
which leads to the MIP formulation of the complete problem given by
\begin{equation}\label{formulsimple2}
 \min\set{y_1\,:\, \eqref{formulsimpleminus},\: \eqref{formulsimpletwo}}.
 \end{equation}
We can check that the optimal value of the LP relaxation of \eqref{formulsimple2} is equal to one. In contrast, we can also check that the optimal value of the LP relaxation of \eqref{formulsimple} is zero. That is, we have constructed a stronger MIP formulation for minimizing a nonlinear objective over a union of polyhedra by directly including the linearization of the objective in our construction procedure.

Given that incorporating additional structure in the ground set can allow us to construct stronger formulations, it seems at first that the optimal approach will be to simply add all constraints. However, this can quickly lead to embedding objects $\En(Q,D,C)$ that are very complex or difficult to study; if $Q$ is restricted to some minimal ``interesting'' substructure, we will see that we are better equipped to study and construct strong formulations.

\section{MIP formulations and valid inequalities for pairwise layouts}\label{sec:twobox}

\subsection{Unary formulation}
We start by analyzing a simple, yet nontrivial, substructure for which we are able to construct a strong (i.e. ideal) formulation. Take $Q^{lb}_{i,j} \defeq \left\{(c_i,c_j,\ell_i,\ell_j) \in \bbR^8: \eqref{eqn:sitb},\;\eqref{eqn:lowerbound} \right\}$;
that is, the set that imposes that the boxes lie completely on the floor and lower bounds on the box widths. Using this set of common constraints, disjunction $D^4$ and the unary encoding  we can construct the following small ideal formulation. Throughout, we will use the notation $\{p,q\}=\{i,j\}$ as enumeration over the two orderings $(i,j)$ and $(j,i)$.
\begin{theorem} \label{thm:unary}
	The following is a formulation for  $\En(Q^{lb}, D^4, U^4)$:
\begin{subequations} \label{eqn:unary-formulation}
\begin{alignat}{2}
    \frac{1}{2}\ell^s_p + lb^s_q u^s_{q,p} \leq c^s_p \leq L^s - \frac{1}{2}\ell^s_p - lb^s_q u^s_{p,q} \quad &\forall s \in \{x,y\}, \{p,q\} = \{i,j\} \label{eqn:tight-sitb} \\
    c^s_p + \frac{1}{2}\ell^s_p \leq c^s_q - \frac{1}{2}\ell^s_q + L^s(1-u^s_{p,q}) \quad &\forall s \in \{x,y\}, \{p,q\}=\{i,j\} \label{eqn:unary-nonoverlap} \\
	\ell^s_p \geq lb^s_p \quad &\forall s \in \{x,y\}, p \in \{i,j\} \label{eqn:unary-last-shared} \\
	u^x_{i,j} + u^x_{j,i} + u^y_{i,j} + u^y_{j,i} = 1 \quad& \label{eqn:sum-to-one} \\
	u^s_{p,q} \in \{0,1\} \quad &\forall s \in \{x,y\}, \{p,q\}=\{i,j\}. \label{eqn:unary-lowerbounds}
\end{alignat}
\end{subequations}
	If $lb^s_i + lb^s_j < L^s$ for both $s \in \{x,y\}$, then this formulation is ideal.
\end{theorem}

We dub \eqref{eqn:unary-formulation} the \emph{unary formulation}. However, we do not use the same naming convention for the binary variables as in formulation \eqref{eqn:unary-bigM}.
Instead we rename $v_1$, $v_2$, $v_3$ and $v_4$ to $u^x_{i,j}$, $u^x_{j,i}$, $u^y_{i,j}$ and $u^y_{j,i}$. The reason for this is that the $0/1$ variables from \eqref{eqn:unary-formulation} have the nice interpretation that $u^s_{i,j} = 1 \Longrightarrow \scrB_i \leftarrow_s \scrB_j$. This interpretation forms the basis for the FLP2 formulation in \cite{Meller:1999} and \cite{Sherali:2003}. In fact, the unary formulation \eqref{eqn:unary-formulation} is very similar to FLP2, but with the addition of the tightened stay-on-the-floor constraints \eqref{eqn:tight-sitb}. In the sequel we use similar naming conventions for the $0/1$ variables when they have helpful interpretations.

\subsection{Binary formulations}\label{sec:binaryform}
The unary encoding uses codes of four bits to differentiate between four choices. If we instead use a binary encoding, we  only need two bits (i.e. codes of length two) to impose this same decision. In contrast to unary encodings, the specific assignment of codes to branches for binary encodings can result in significantly different formulations \citep{Vielma:2015a}. However, because of symmetry, for binary encodings of length two we may restrict our attention to two possible choices. The first encoding corresponds to the unique (up to symmetry) Gray code \citep{Savage:1997} with two bits given by  $GB^4 \defeq \left\{ (0,0), (1,0), (0,1), (1,1) \right\}$,
and the second corresponds to the codes $BB^4 \defeq \left\{ (0,0), (1,1), (1,0), (0,1) \right\}$.
Both choice of codes and their corresponding encodings can be used to reinterpret existing formulations from the literature. The following proposition shows that the big-$M$ approach can be used to construct a simple formulation for the Gray encoding, which can be seen as the basis of formulation FLP-SP introduced in \cite{Meller:2007}.

\begin{proposition}\label{graybigm}
	A valid formulation for $\En(Q^{lb},D^4,GB^4)$ is:
\begin{subequations} \label{eqn:gray-binary-formulation}
\begin{alignat}{3}
	\frac{1}{2}\ell^s_p \leq c^s_p &\leq L^s - \frac{1}{2}\ell^s_p &\quad& \forall s \in \{x,y\}, p \in \{i,j\} \label{eqn:gb-first} \\
    \ell^s_p &\geq lb^s_p &\quad& \forall s \in \{x,y\}, p \in \{i,j\} \\
	c^y_i + \frac{1}{2}\ell^y_i &\leq c^y_j - \frac{1}{2}\ell^y_j + L^y(w_1+w_2) \\
	c^x_i + \frac{1}{2}\ell^x_i &\leq c^x_j - \frac{1}{2}\ell^x_j + L^x(1-w_1+w_2) \\
	c^y_j + \frac{1}{2}\ell^y_j &\leq c^y_i - \frac{1}{2}\ell^y_i + L^y(2-w_1-w_2) \\
	c^x_j + \frac{1}{2}\ell^x_j &\leq c^x_i - \frac{1}{2}\ell^x_i + L^x(1+w_1-w_2) \\
	w &\in \{0,1\}^2 \label{eqn:gb-last}
\end{alignat}
\end{subequations}
\end{proposition}

Formulation FLP-SP is obtained from \eqref{eqn:gray-binary-formulation} by adding the ``sequence-pair'' inequalities for the $N$-box formulation introduced in \citep{Meller:2007}. For completeness, we present the sequence pair inequalities in Appendix \ref{app:sequence-pair}. We will see in Section \ref{sec:computations} that the FLP-SP is the most competitive formulation from the literature on our computational benchmarks.

If instead we attempt to construct a formulation for $\En(Q^{lb},D^4,BB^4)$, we can easily reconstruct the BLDP1 formulation from \cite{Castillo:2005}, which we present in Appendix \ref{app:bldp1}.

\subsection{Refined disjunction formulation} \label{sec:refined-disjunction}
\begin{figure}[htpb]
	\centering
	\includegraphics[width=.45\linewidth]{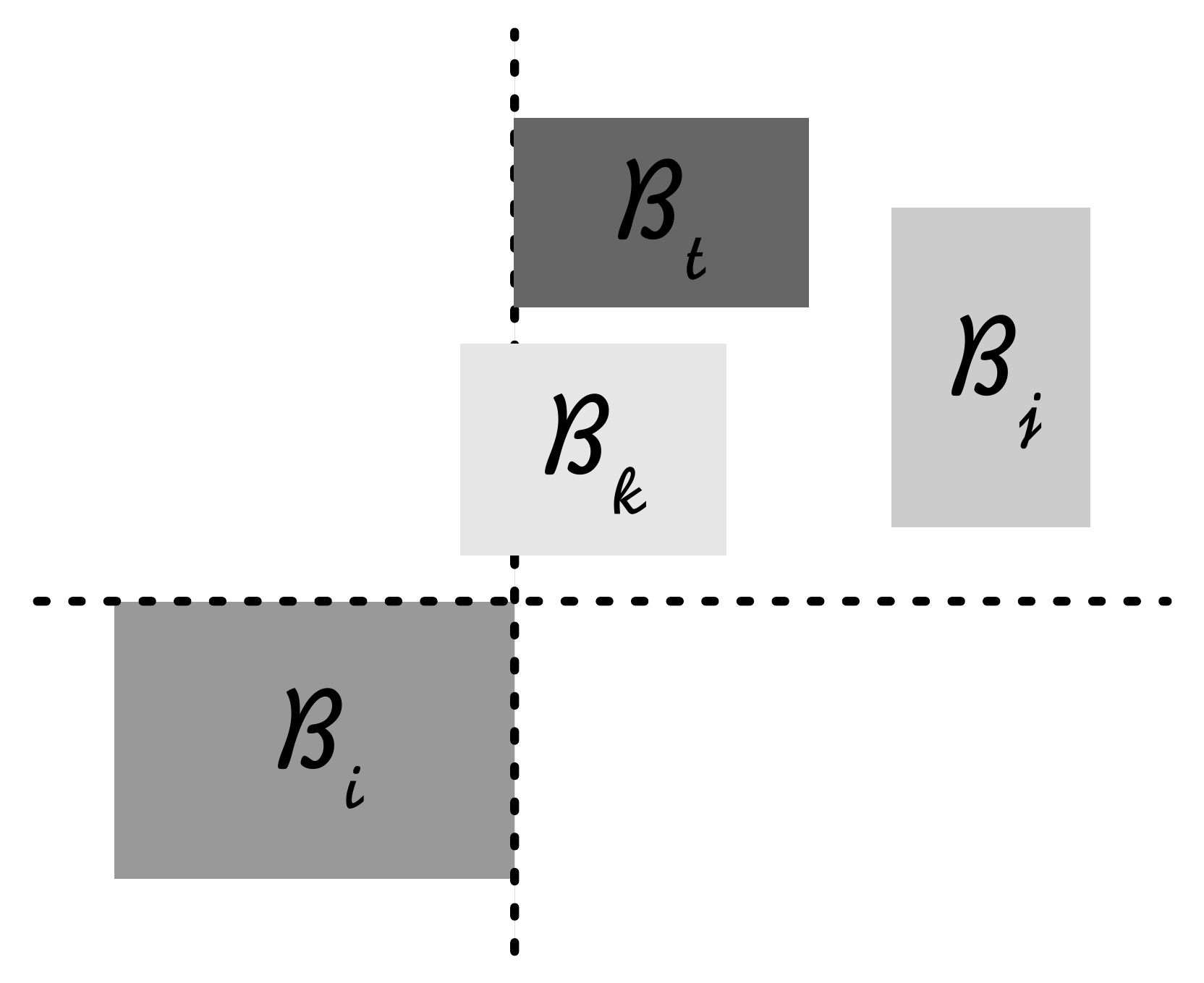}\quad
	\includegraphics[width=.45\linewidth]{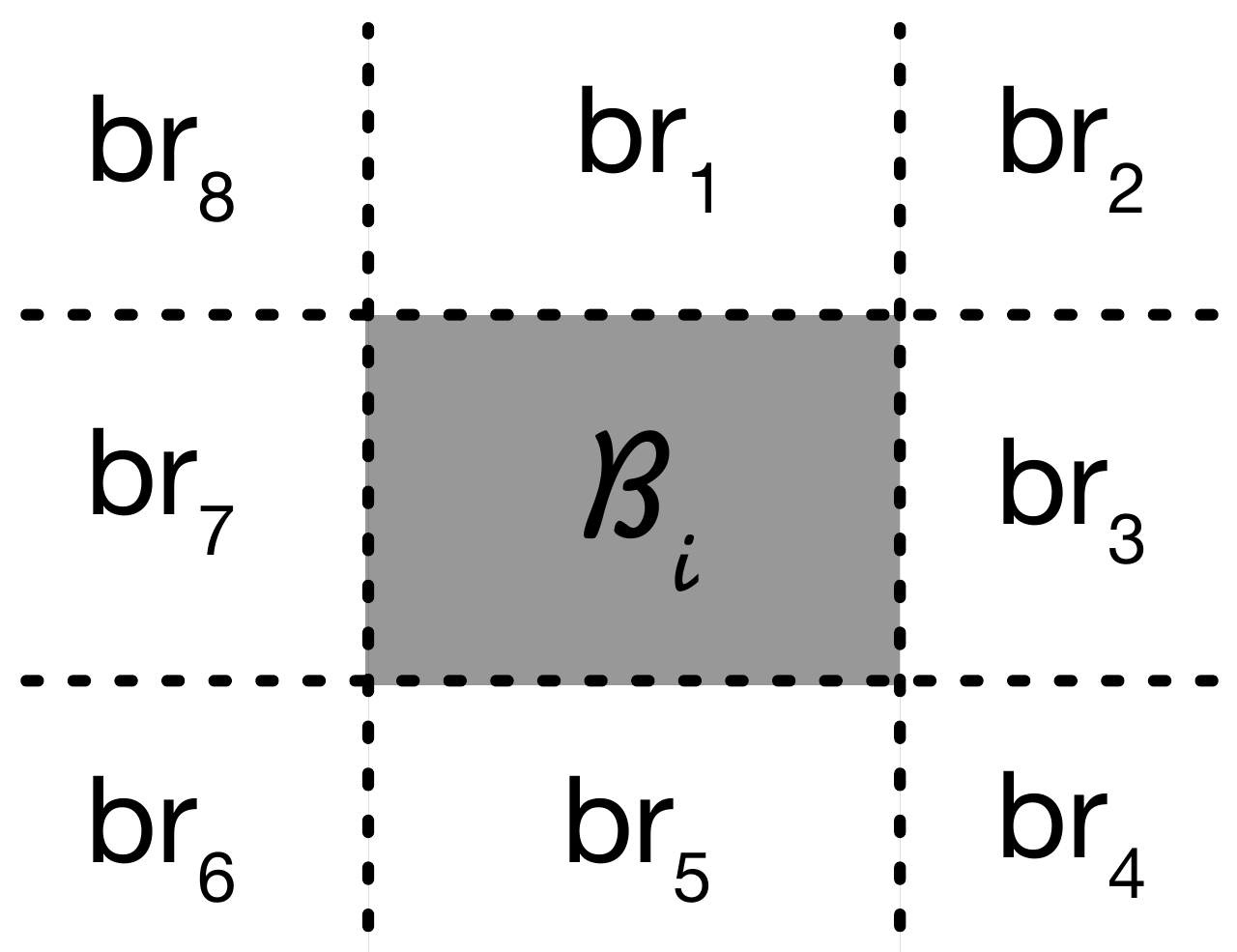}
	\caption{(Left) An illustration of the possible redundancies. In the depiction, $\scrB_i$ precedes $\scrB_j$, $\scrB_k$, and $\scrB_t$ in direction $y$. It also precedes $\scrB_j$ in direction $x$, does not precede $\scrB_k$ in direction $x$, and both precedes and does not precede $\scrB_t$ in direction $x$. \\(Right) The eight branches of the disjunction $D^8$, illustrated via the relative position of $\scrB_j$ to $\scrB_i$.}
	\label{fig:8-configurations}
\end{figure}

While the disjunction $D^4_{i,j}$ is sufficient to enforce that $\scrB_i$ and $\scrB_j$ do not overlap, its simplicity has a downside when used in a MIP framework. The disjunction is not sufficiently refined in the sense that there exist many feasible layouts that satisfy multiple branches at once. For example, in Figure \ref{fig:8-configurations}, we see that $\scrB_i$ precedes $\scrB_j$ in both the $x$ and $y$ directions. Therefore, in any embedding constructed using $D^4_{i,j}$, there exist two points that project down to the same layout (that is, they differ only in their assigned codes). In practice, this redundancy can hamper the progress of branch-and-bound solvers, which must explicitly enumerate these solutions (and all nodes preceding them in the tree) to prove optimality.

To help remove this redundancy from the feasible set, we present a refined disjunction that is logically equivalent to $D^4$. In Definition \ref{defn:precede}, we presented a linear inequality that enforces that $\scrB_i$ precedes $\scrB_j$. For our refined disjunction, we will need a description of the opposite.
\begin{definition}
We say that $\scrB_i$ \emph{does not precede} $\scrB_j$ (denoted by $\scrB_i \not\leftarrow_s \scrB_j$) if $c^s_i + \frac{1}{2}\ell^s_i \geq c^s_j - \frac{1}{2}\ell^s_j$.
\end{definition}
Referring back to Figure \ref{fig:8-configurations}, we see that $\scrB_i$ precedes $\scrB_k$ in direction $y$, but does not precede $\scrB_k$ in direction $x$ (and vice versa). Note in particular that, if $c^s_i + \frac{1}{2}\ell^s_i = c^s_t - \frac{1}{2}\ell^s_t$ (as with $\scrB_i$ and $\scrB_t$ in Figure \ref{fig:8-configurations}), we have that both $\scrB_i \leftarrow_s \scrB_j$ and $\scrB_i \not\leftarrow_s \scrB_j$ simultaneously.


With the two definitions, we can construct a refinement of $D^4_{i,j}$ given by
$D^8_{i,j} \defeq \bigvee_{k=1}^8 br^k_{i,j}$ where
\begin{alignat*}{2}
	br^1_{i,j} &= (\scrB_i \leftarrow_y \scrB_j) \wedge (\scrB_i \not\leftarrow_x \scrB_j) \wedge (\scrB_j \not\leftarrow_x \scrB_i), \quad\quad
	br^2_{i,j} &= (\scrB_i \leftarrow_y \scrB_j) \wedge (\scrB_i \leftarrow_x \scrB_j) \\
	br^3_{i,j} &= (\scrB_i \leftarrow_x \scrB_j) \wedge (\scrB_i \not\leftarrow_y \scrB_j) \wedge (\scrB_j \not\leftarrow_y \scrB_i), \quad\quad
	br^4_{i,j} &= (\scrB_i \leftarrow_x \scrB_j) \wedge (\scrB_j \leftarrow_y \scrB_i) \\
	br^5_{i,j} &= (\scrB_j \leftarrow_y \scrB_i) \wedge (\scrB_i \not\leftarrow_x \scrB_j) \wedge (\scrB_j \not\leftarrow_x \scrB_i), \quad\quad
	br^6_{i,j} &= (\scrB_j \leftarrow_x \scrB_i) \wedge (\scrB_j \leftarrow_y \scrB_i) \\
	br^7_{i,j} &= (\scrB_j \leftarrow_x \scrB_i) \wedge (\scrB_i \not\leftarrow_y \scrB_j) \wedge (\scrB_j \not\leftarrow_y \scrB_i), \quad\quad
	br^8_{i,j} &= (\scrB_j \leftarrow_x \scrB_i) \wedge (\scrB_i \leftarrow_y \scrB_j).
\end{alignat*}
We have taken a refinement of $D^4$ by splitting the regions satisfying two branches at once into the new branches $br^2,br^4,br^6$, and $br^8$, and shrinking the other branches to exclude these new regions. See Figure \ref{fig:8-configurations} for an illustration.

With 8 branches in the disjunction we need  codes of length at least $3=\log_2(8)$. However, in lieu of chasing the formulation with the smallest number of $0/1$ variables (i.e. a binary formulation), we instead take the encoding $C^8 \defeq \left\{ \e^1, \e^1+\e^2, \e^2, \e^2+\e^3, \e^3, \e^3+\e^4, \e^4, \e^4+\e^1 \right\} \subseteq \{0,1\}^4$.
Intuitively, we have taken the codes from the unary encoding for the regions $br^1,br^3,br^5$, and $br^7$, and taken the codes for the new branches as the sum of the codes assigned to the two branches in $D^4$ the region satisfies. For example, we get $br^2$ by taking the intersection of $d^1$ (code $\e^1$) and $d^3$ (code $\e^2$), so we take the corresponding code as $\e^1+\e^2$. We see in the following proposition that by ``shadowing'' the unary embedding in this way, we are able to construct a formulation for $\En(Q^{lb},D^8,C^8)$ that is very similar to the ideal formulation \eqref{eqn:unary-formulation} for $\En(Q^{lb},D^4,U^4)$.

\begin{proposition} \label{prop:refined}
	The following is a valid formulation for $\En(Q^{lb},D^8,C^8)$:
    \begin{subequations} \label{eqn:refined-unary-formulation}
	\begin{align}
	    \frac{1}{2}\ell^s_p + lb^s_q z^s_{q,p} \leq c^s_p \leq L^s - \frac{1}{2}\ell^s_p - lb^s_q z^s_{p,q} \quad &\forall s \in \{x,y\}, \{p,q\} = \{i,j\} \label{eqn:ru-tight-sitb} \\
        c^s_p + \frac{1}{2}\ell^s_p \leq c^s_q - \frac{1}{2}\ell^s_q + L^s(1-z^s_{p,q}) \quad &\forall s \in \{x,y\}, \{p,q\}=\{i,j\} \label{eqn:ru-nonoverlap} \\
	    \ell^s_p \geq lb^s_p \quad &\forall s \in \{x,y\}, p \in \{i,j\} \label{eqn:ru-lowerbounds} \\
		z^x_{i,j} + z^x_{j,i} + z^y_{i,j} + z^y_{j,i} \geq 1 \quad & \\
		z^s_{i,j} + z^s_{j,i} \leq 1 \quad &\forall s \in \{x,y\} \label{eqn:refined-proof-helper} \\
		z^s_{p,q} \in \{0,1\} \quad &\forall s \in \{x,y\}, \{p,q\}=\{i,j\} \label{eqn:refined-keeper-last} \\
		c^s_p + \frac{1}{2}\ell^s_p + L^sz^s_{p,q} \geq c^s_q - \frac{1}{2}\ell^s_q + (lb^s_p+lb^s_q)(z^s_{i,j} + z^s_{j,i}) \quad &\forall s \in \{x,y\}, \{p,q\}=\{i,j\} \label{eqn:refined-new}.
	\end{align}
    \end{subequations}
\end{proposition}

We dub formulation  \eqref{eqn:refined-unary-formulation} the \emph{refined unary formulation}, and conjecture that it is the strongest possible for $\En(Q^{lb},D^8,C^8)$.

\begin{conjecture}
	Formulation \eqref{eqn:refined-unary-formulation} is ideal.
\end{conjecture}

Finally, we note that our choice of codes induce the following nice interpretation for the $0/1$ variables:
\begin{subequations} \label{eqn:z-variable-implication}
\begin{align}
    z^s_{i,j} &= 0 \Longrightarrow \scrB_i \not\leftarrow_s \scrB_j \label{eqn:z-variable-implication-zero} \\
	z^s_{i,j} &= 1 \Longrightarrow \scrB_i \leftarrow_s \scrB_j.
\end{align}
\end{subequations}
We note that this is a stronger interpretation than is possible for the unary formulation, for which the implication \eqref{eqn:z-variable-implication-zero} is not necessarily true.


\section{Constructing valid inequalities for embeddings}\label{sec:inequalities}
In Section \ref{sec:formulations} we have seen how the embedding approach can be used to construct valid formulations for substructures of the pairwise FLP. In particular, we chose a subset of variables and constraints ($Q^{lb}$) for which the analysis is tractable. However, in Section \ref{sec:common-constraints} we have seen that incorporating more of the common constraint structure in the ground set $Q$ can allow us to construct much stronger formulations. Therefore, in this section we explore embeddings of more complex substructures $Q$. However, since the facial structure of the embedding objects grows considerably more complex, we only focus on constructing valid inequalities for these new embeddings.

In the remaining sections, we will express all inequalities for the refined unary encoding. Fortunately, it is sometimes possible to translate valid inequalities between different encodings. We now present a self contained description of such translations for the FLP.

\begin{proposition} \label{prop:map-from-refined}
    Let $Q^{FLP} \defeq \{(c_i,c_j,\ell_i,\ell_j) \in \bbR^{8}: \eqref{eqn:sitb}, \eqref{eqn:area}, \eqref{eqn:bounds}\}$ and consider an inequality $a^Tc + b^T\ell + d^Tz \leq f$ with $d \geq 0$ that is valid for $\En(Q^{FLP},D^8,C^8)$. Then
    \begin{itemize}
        \item $a^Tc + b^T\ell + d^T\scrA^U(u) \leq f$ is valid for $\En(Q^{FLP},D^4,U^4)$, and
	    \item $a^Tc + b^T\ell + d^T\scrA^{GB}(w) \leq f$ is valid for $\En(Q^{FLP},D^4,GB^4)$,
    \end{itemize}
    where $\scrA^U(u) = u$ is the affine mapping that identifies $z^s_{p,q}$ with $u^s_{p,q}$ and
\[
    \scrA^{GB}(w) \defeq \left(-w_1-w_2+1,\:\: w_1-w_2,\:\: w_1+w_2-1,\:\: -w_1+w_2 \right).
\]
is the affine mapping that identifies $(z^y_{i,j},z^x_{i,j},z^y_{j,i},z^x_{j,i})$ with $\scrA^{GB}(w)$.
\end{proposition}

\begin{proposition} \label{prop:map-to-refined}
    If inequality $a^Tc + b^T\ell + d^Tu \leq f$ is valid for $\En(Q^{FLP},D^4,U^4)$, then $a^Tc + b^T\ell + d^Tz \leq f$ is valid for $\En(Q^{FLP},D^8,C^8)$ if either $d^y_{i,j} = d^y_{j,i} = 0$ or $d^x_{i,j} = d^x_{j,i} = 0$.
\end{proposition}

\subsection{Upper bound inequalities}

In the previous section we chose the base set $Q^{lb}$ such that only lower bounds on the widths were included in the formulation. This was to make the formulation analysis tractable, but enforcing the aspect-ratio constraints via \eqref{eqn:bounds} naturally includes upper bounds as well. Therefore, we can consider the set $\En(Q^{ub},D^8,C^8)$ induced by $Q^{ub} \defeq \{(c,\ell) \in Q^{lb}: \eqref{eqn:upperbound}\}$.
\begin{proposition} \label{prop:ub-cuts}
	For any assignments $\{r,s\} = \{x,y\}$ and $\{p,q\} = \{i,j\}$, then
	\begin{equation}
		c^s_p + ub^s_q(1-z^s_{q,p}) \geq \frac{1}{2}\ell^s_p + \ell^s_q \label{eqn:ub-c1} \\
	\end{equation}
	is a valid inequality for $\En(Q^{ub},D^8,C^8)$. If $L^s < ub^s_p + ub^s_q$,
	\begin{equation}
		z^r_{p,q} + z^r_{q,p} \geq \frac{\ell^s_p + \ell^s_q - L^s}{ub^s_p + ub^s_q - L^s} \label{eqn:ub-c2}
	\end{equation}
	is valid for both $\En(Q,D^4,U^4)$ and $\En(Q,D^8,C^8)$.
\end{proposition}

\subsection{Objective inequalities} \label{sec:objective}
The objective \eqref{eqn:objective} is nonlinear but is straightforward to linearize in the usual fashion with auxiliary variables $(d^x_{i,j},d^y_{i,j})$ and the constraints
\begin{equation} \label{eqn:linearized-objective}
	d^s_{i,j} \geq c^s_i - c^s_j, \quad\quad
	d^s_{i,j} \geq c^s_j - c^s_i.
\end{equation}
Even though this type of linearization is a very common MIP formulation technique, it is often not incorporated into polyhedral studies explicitly. To do this for the pairwise FLP, consider the augmented base set $Q^{obj}_{i,j} = \left\{(c_i,c_j,\ell_i,\ell_j,d_{i,j}) \in \bbR^{4+4+2}: \eqref{eqn:sitb},\eqref{eqn:lowerbound},\eqref{eqn:linearized-objective}\right\}$.
The resulting encoding $\En(Q^{obj},D^8,C^8)$ leads to a collection of inequalities that serve to lower bound the auxiliary objective variables $d$.
\begin{proposition} \label{prop:obj-cuts}
	Choose $s \in \{x,y\}$ and some assignment $\{p,q\} = \{i,j\}$. Then the following are valid inequalities for $\En(Q^{obj},D^8,C^8)$:
	\begin{alignat}{3}
		d^s_{i,j} &\geq \frac{1}{2}(\ell^s_i+\ell^s_i) - L^s(1-z^s_{i,j}-z^s_{j,i}) \label{eqn:obj1} \\
		d^s_{i,j} &\geq c^s_p - c^s_q + \ell^s_p + lb^s_q(z^s_{p,q} + z^s_{q,p}) - L^s(1-z^s_{p,q}) \label{eqn:obj2} \\
		d^s_{i,j} &\geq c^s_p - c^s_q + (lb^s_p+lb^s_q)z^s_{p,q} \label{eqn:obj3} \\
		2d^s_{i,j} &\geq \ell^s_p - L^s(1-z^s_{p,q}-z^s_{q,p}) + lb^s_q(z^s_{p,q}+z^s_{q,p}) \label{eqn:obj4}
	\end{alignat}
\end{proposition}

Note that we are now adding both constraints and variables to our ground set $Q^{obj}$. These inequalities are especially significant, since they explicitly incorporate the objective function, and the MIP relaxation lower bounds for the FLP are quite poor (see Section \ref{sec:lowerbounds}).

\section{From pairwise to $N$ boxes}\label{sec:multibox}
Thus far we have only considered representations for $\hat{\scrL}_{i,j}$, the relationships between a single pair of boxes. In this section we address how to use the results derived for the pairwise formulations to construct strong formulations for the original $N$-box floor layout problem.

\subsection{Multi-box formulations} \label{sec:multibox-formulation}

Since all the constraints for the FLP involve at most two boxes, it suffices to consider each pair of boxes separately, construct a pairwise formulation, and identify all repeated variables across these pairwise formulations as follows.

\begin{proposition} \label{prop:pairwise-to-Nbox}
	Consider pairwise formulations $F^{i,j}$ for each pair of boxes $(i,j) \in \scrP$ over the variables $(c_i,c_j,\ell_i,\ell_j,v^{i,j}) \in \bbR^{8} \times \{0,1\}^{m_{i,j}}$. If $M \defeq \sum_{(i,j) \in \scrP} m_{i,j}$, then $\left\{(c,\ell,v) \in \bbR^{4N} \times \{0,1\}^M: (c_i,c_j,\ell_i,\ell_j,v^{i,j}) \in F^{i,j} \quad \forall (i,j) \in \scrP \right\}$ is a formulation for $\scrL$.
\end{proposition}

In particular, if we take the refined unary formulation for each pair of boxes, we construct the following formulation for the $N$-box FLP.
\begin{corollary}\label{prop:pairwise-to-NboxC}
	Take $F^{RU}_{i,j} = \{(c_i,c_j,\ell_i,\ell_j,v^{i,j}) \in \bbR^{8} \times \{0,1\}^4: \eqref{eqn:refined-unary-formulation}\}$. Then $F^{RU} \defeq \left\{ (c,\ell,v) \in \bbR^{8N} \times \{0,1\}^{2N(N-1)} : (c_i,c_j,\ell_i,\ell_j,v^{i,j}) \in F^{RU}_{i,j} \: \forall (i,j) \in \scrP \right\}$ is a valid formulation for $\scrL$.
\end{corollary}

While this approach is sufficient to construct a valid formulation, constructing  disjunctive and MIP formulation for multiple pairs of boxes can lead to stronger formulations. However, such formulations can be significantly larger and/or more complicated. For this reason we instead concentrate on identifying valid inequalities for such multi-pair or multi-box formulations to strengthen the single-pair formulation from Propositions~\ref{prop:pairwise-to-Nbox} and Corollary~\ref{prop:pairwise-to-NboxC}.
\subsection{Multi-box cutting planes} \label{ss:multibox-cuts}
When working with more than two boxes at once, the notion of spatial transitivity appears; that is, for any $s\in \set{x,y}$ we have $\scrB_i \rightarrow_s \scrB_t \rightarrow_s \scrB_j \Longrightarrow \scrB_i \rightarrow_s \scrB_j$.
We can use this property to generalize many of the pairwise valid inequalities introduced thus far to the multi-box setting, in a similar way to Section 3 of \cite{Meller:1999}.

Consider a pair of boxes $(i,j) \in \scrP$ and take an arbitrary path $P = \{(t^0,t^1),\ldots,(t^m,t^{m+1})\} \subseteq \scrP$,
where $t^0 = i$ and $t^{m+1} = j$. We define an affine function of the form $\scrM^s_P(z) \defeq 1 + \sum_{\xi=1}^{m+1} \left(z^s_{t^{\xi-1},t^\xi} - 1\right)$.
Our function enjoys the following property:
\begin{equation}
	\scrM^s_P(z) \begin{cases} \label{eqn:multibox-map-property}
		= 1 & \scrB_{t^0} \leftarrow_s \scrB_{t^1} \leftarrow_s \cdots \leftarrow_s \scrB_{t^{m+1}} \\
		\leq 0 & \text{otherwise}.
	\end{cases}
\end{equation}
This will function as an (underestimator for the) indicator function for when we have a particular chain of boxes $P$ along direction $s$. We can use this to extend the logic of the pairwise inequalities we have developed. For a simple example, if $\scrB_i \leftarrow_s \scrB_t \leftarrow_s \scrB_j$, then we know that $\scrB_i$ and $\scrB_j$ are separated in direction $s$ by at least the smallest width $\scrB_t$ can take along that direction, and so $c^s_i + \frac{1}{2}\ell^s_i + lb^s_t \leq c^s_j - \frac{1}{2}\ell^s_j$.
This tightening can be exploited in the inequalities derived previously, leading a host of new valid inequalities for the multi-box FLP.
\begin{proposition} \label{prop:multibox-cuts}
	Consider the pair $(\scrB_i,\scrB_j)$ and an arbitrary path $P = \{(t^0,t^1),\ldots,(t^m,t^{m+1})\}$, where $i=t^0$ and $j=t^{m+1}$ and $m \geq 1$. Choose assignments $\{r,s\} = \{x,y\}$ and $\{p,q\} = \{i,j\}$ and define $\gamma_P \defeq \sum_{\xi=1}^m lb^s_{t^\xi}$. Then the following are valid inequalities for $F^{RU}$:
\begin{alignat}{3}
    d_{i,j}^s &\geq \frac{1}{2}(\ell_i^s + \ell_j^s) - L^s(1-z_{i,j}^s-z_{j,i}^s) + \gamma_P\scrM^s_P(z) \label{eqn:multi1} \\
    d_{i,j}^s &\geq c_i^s - c_j^s + \ell_p^s + lb_q^s(z_{i,j}^s+z_{j,i}^s) - L^s(1-z_{p,q}) + \gamma_P\scrM^s_P(z) \label{eqn:multi2} \\
    d^s_{i,j} &\geq c^s_i - c^s_j + (lb^s_i+lb^s_j)z^s_{i,j} + \gamma_P\scrM^s_P(z) \label{eqn:multi3} \\
    2d_{i,j}^s &\geq \ell_p + lb_q^s(z_{i,j}^s+z_{j,i}^s) - L^s(1-z_{i,j}^s-z_{j,i}^s) + 2\gamma_P\scrM^s_P(z) \label{eqn:multi4} \\
    \frac{1}{2}\ell_j^s + lb^s_iz^s_{i,j} + \gamma_P\scrM^s_P(z) &\leq c^s_j \label{eqn:multi5} \\
    c_i^s + \gamma_P\scrM^s_P(z) &\leq L^s - \frac{1}{2}\ell^s_i - lb^s_jz^s_{i,j} \label{eqn:multi6} \\
    c^s_i + \frac{1}{2}\ell^s_i + \gamma_P\scrM^s_P(z) &\leq c_j^s - \frac{1}{2}\ell_j^s + L^s(1 - z_{i,j}^s) \label{eqn:multi7}.
 \end{alignat}

\end{proposition}
Proposition \ref{prop:multibox-cuts} provides an exponential number of valid inequalities for the $N$-box FLP. For small paths (e.g. $|P| = 2$), these inequalities can be added to the formulation directly; this is the approach we take in the computational trials.

\section{Computational results} \label{sec:computations}
In the computational trials we compare four formulations, along with four different collections of valid inequalities added to each formulation. The unary formulation, denoted \texttt{U}, is based on the pairwise unary formulation \eqref{eqn:unary-formulation}; this is a strengthened version of the FLP2 formulation from \cite{Meller:1999}. The \texttt{BLDP1} formulation from \cite{Castillo:2005a} is also tested; see Appendix~\ref{app:bldp1}. The third formulation is the sequence-pair formulation (\texttt{SP}) from \cite{Meller:2007}, derived by adding global constraints to a formulation derived from \eqref{eqn:gray-binary-formulation}; see Appendix~\ref{app:sequence-pair}. Finally, we compare with the new refined unary formulation \eqref{eqn:refined-unary-formulation}, which we denote \texttt{RU}. We note that we observe a slight computational advantage for using the simple stay-on-the-floor constraints \eqref{eqn:sitb} rather than the tightened versions \eqref{eqn:tight-sitb} and \eqref{eqn:ru-tight-sitb}, since we may aggregate them and add a single copy, rather than one for each pair $(i,j) \in \scrP$; we instead add these tightened constraints as valid inequalities.

We will compare each of these formulations with one of four levels of valid inequalities added to the formulation (that is, they are not separated dynamically). The first will be no valid inequalities. The second will use the V2 and B2 families of inequalities appearing in \cite{Meller:1999}; the formulation name will be appended with \texttt{+} if these inequalities are added. We present these inequalities in Appendix \ref{app:literature-cuts} for completeness. The \texttt{VI} tag will be used for formulations with the new inequalities derived in this work added. In particular, we use (\ref{eqn:obj1}-\ref{eqn:obj4}), \eqref{eqn:ub-c1}, and (\ref{eqn:multi5}-\ref{eqn:multi7}) for paths $|P|=2$. For the \texttt{RU} formulation, we also add \eqref{eqn:ub-c2} and \eqref{eqn:ru-tight-sitb}. Adding all of these to the formulation proved impractical in the branch-and-bound setting, so we instead add an $O(n)$ subset of these inequalities with \texttt{VI}\footnote{Specifically, we add the $N$ pairs $(i,j)$ with largest objective coefficient $p_{i,j}$. For the three-box inequalities, we choose the $N$ triplets $(i,j,k)$ which maximize $p_{i,j} + p_{i,k} + p_{j,k}$, and add the inequalities corresponding to all 6 paths (i.e. permutations) through $(i,j,k)$}. Finally, we consider adding (an $O(n)$ subset of) the ``multi-box objective cuts'' (\ref{eqn:multi1}-\ref{eqn:multi4}) for paths $P$ of length 3, and denote this by appending a \texttt{3}.

There is a standard symmetry-breaking approach presented in \cite{Sherali:2003} that we will use in all the computational examples to follow (except for the relaxation gap discussion in Section \ref{sec:lowerbounds}, where we will discuss the effect of the symmetry-breaking explicitly). We present the symmetry-breaking scheme in Appendix \ref{app:symmetry-breaking} for completeness.

For our benchmarks, we use the \texttt{hp} (11 boxes), \texttt{apte} (9 boxes), and \texttt{xerox} (10 boxes) benchmarks from the MCNC benchmark collection \cite{North-Carolina:2015}. Additionally, we will use the \texttt{Armour62-1} (20 boxes) and \texttt{Armour62-2} (20 boxes) instances from \cite{Armour:1962}, the \texttt{Bazaraa75-1} (13 boxes) and \texttt{Bazaraa75-2} (14 boxes) instances from \cite{Bazaraa:1975}, the \texttt{Camp91} (10 boxes) instance from \cite{Camp:1991}, the \texttt{Bozer91} (15 boxes) instance from \cite{Bozer:1991}, and the \texttt{Bozer97-1} (9 boxes) and \texttt{Bozer97-2} (12 boxes) instances from \cite{Bozer:1997} instances. To the best of our knowledge, none of these instances have been solved to optimality before in the literature. From each of these 11 base instances, we create a family of related instances by 1) selecting the aspect ratio $\alpha \in \{4,5,6\}$, and adding three possible levels random noise to the nonzero problem data (perturbations of the form $x \leftarrow (1 + \gamma t) x$ for standard normal $t$ and for $\gamma \in \{0.0,0.1,0.2\}$). For the remainder, we refer to each instance according to the schema \texttt{instance\_name-$\gamma$($\alpha$)}; for example, the \text{xerox} instance with aspect ratio $\alpha=4$ and perturbation factor $\gamma=0.1$ is \texttt{xerox-0.1(4)}.

To construct the formulations and interface with the solver, we use the JuMP algebraic modeling language from \cite{Dunning:2015a,Lubin:2013}; JuMP is written in the Julia programming language (see \cite{Bezanson:2012}). We performed the experiments on an Intel i7-3770 3.40GHz Linux workstation with 32GB of RAM. All trials use CPLEX v12.6 with a maximum runtime of 4 hours. We also performed trials with Gurobi v6.0, but the performance was not competitive. We force CPLEX to use the linearization of the second-order cone constraints ($\texttt{CPX\_PARAM\_MIQCPSTRAT}=2$), as solving the nonlinear problem at the nodes was far slower.

The code used for these computational studies, as well as the benchmark instances in MPS format are available at \url{https://github.com/joehuchette/floor-layout}.

\subsection{Relaxation bound} \label{sec:lowerbounds}
First, we compare the lower bound produced by solving the continuous relaxation of the formulations, with and without valid inequalities added. We present the relative gap percentage $100\frac{U-L}{U}$, where $L$ is the given relaxation lower bound and $U$ is the cost of the best available feasible solution. Note in particular that a relative gap percentage of $100\%$ implies that the relaxation bound of $0$, which is the worst possible for any of the formulations presented here for the FLP. We observe that the formulations naturally fall into two groups with respect to the quality of their relaxation bound, and we summarize the results below (we include a table in Appendix \ref{app:relaxation-gap-table} for completeness).

First, we observe that the ``two-bit'' formulations (\texttt{BLDP1} and \texttt{SP}) have a relaxation gap of $100\%$, even with \emph{all inequalities} discussed in the previous subsection added to the formulation (\texttt{+VI3}). With the symmetry-breaking constraints added to the formulation, the relaxation lower bound is no longer zero, and so the relaxation gap improves slightly (mean 89.3\%, with standard deviation $5.85\%$). The ``four-bit'' formulations (\texttt{U} and \texttt{RU}) also produce a trivial relaxation gap of $100\%$, but adding the valid inequalities helps improve the lower bound considerably (mean $57.4\%$ with standard deviation $8.4\%$). In particular, we can isolate the B2 inequalities from \cite{Meller:1999} and \eqref{eqn:obj3} as the crucial additions to the improvement in the gap.

We also compare the relative gap attained at the root node (with respect to the best known feasible solution) after CPLEX is able to apply advanced techniques such as general purpose cuts and preprocessing.  This improves the gap by roughly 1\%-5\% for most trials. However, there is still an appreciable difference in gap between the ``four bit'' formulations \texttt{U} and \texttt{RU} and the ``two bit'' formulations \texttt{SP} and \texttt{BLDP1}. A complete table is available in Appendix~\ref{app:relaxation-gap-table}.


\subsection{Branching behavior}
\begin{figure}[htpb]
    \centering
    \includegraphics[width=.8\linewidth]{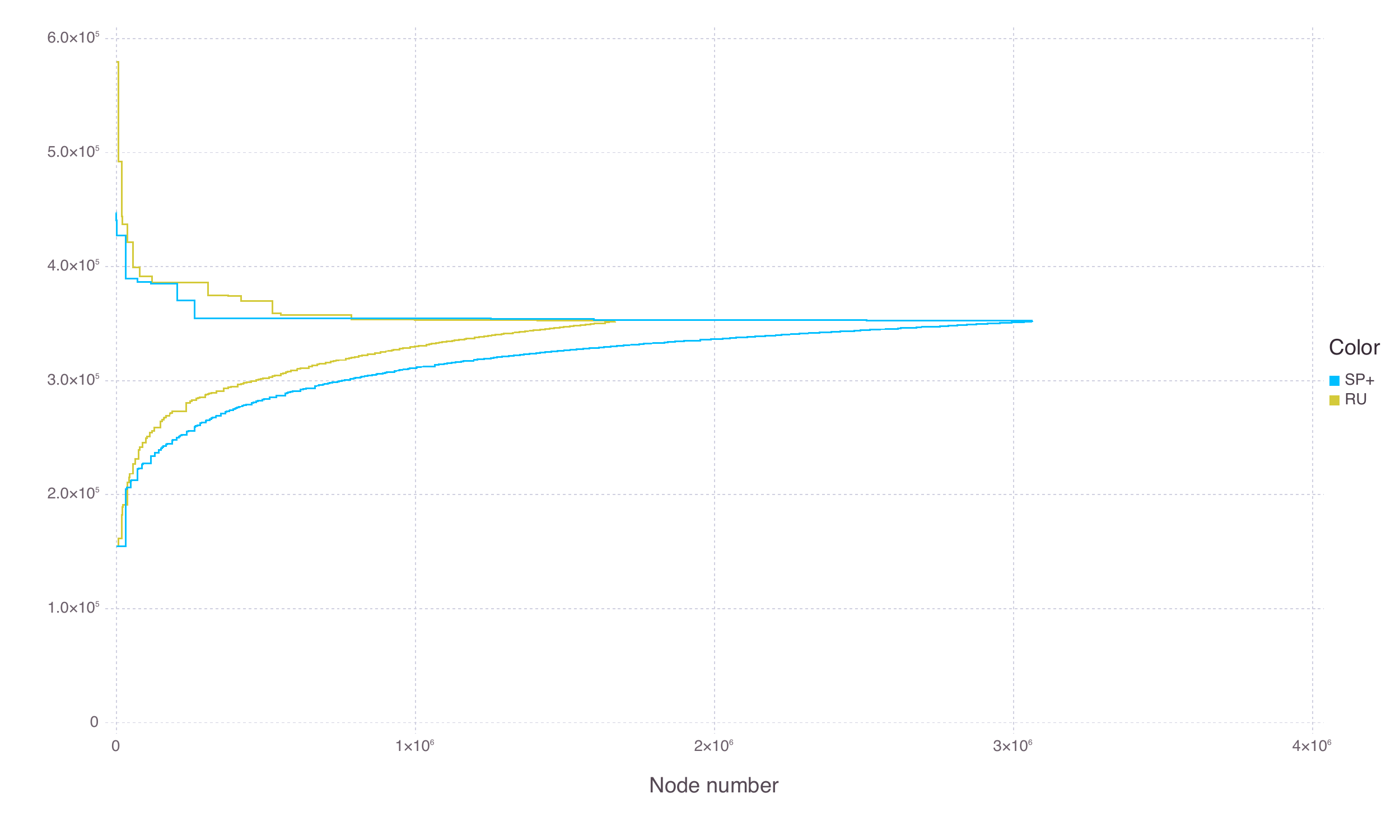} 
    \caption{Plots of lower and upper bounds on optimal cost for the sequence pair (\texttt{SP+}) and refined unary (\texttt{RU}) formulations as a function of node count for the \texttt{xerox-0.0(5)} benchmark instance.}
    \label{fig:branching}
\end{figure}

The rationale for introducing the refined unary formulation \eqref{eqn:refined-unary-formulation} was that many feasible layouts will have multiple corresponding points in the encoding constructed using $D^4$. The refined partition removes many of these redundant solutions, which helps in the branch-and-bound setting, much in the same way symmetry breaking removes equivalent feasible solutions that would otherwise have to be explicitly enumerated in the optimization procedure. Qualitatively, we observe this change of behavior in Figure \ref{fig:branching}, where we compare the progress of the \texttt{SP+} and \texttt{RU} formulations as a function of node count on the \texttt{xerox} benchmark. That is, we compare both the upper and lower bound for both formulations; when they are equal, the solver has proven optimality. We see that the \texttt{RU} formulation requires fewer nodes to prove optimality, as expected. More broadly, this illustration shows the typical trajectory when solving an instance of the FLP: finding a good (often near-optimal) feasible solution early in the procedure, and then steadily improving the lower bound which is far from optimal, until the gap is finally closed and optimality is proven.

\subsection{Solution time}
While the \texttt{RU} formulation offers advantages when looking at the progress with respect to the node count, the advantage is not so clear-cut when looking at solution time. Tables~\ref{table:solve-time-1} and \ref{table:solve-time-2} shows the solution time for each of our benchmarks and formulation and inequality combination. On some harder instances, some approaches overflow the memory on our test machine and terminate prematurely. We denote these runs with a \texttt{ME} for ``Memory Error'', and note that at termination, none of these runs produced bounds that were competitive with the other methods for that particular benchmark instance.

We observe that, in a majority of instances (55 of 99), the new techniques presented in the paper yield the best performing approach. We quickly see that the \texttt{U} and \texttt{U+} formulations are never competitive on the benchmark set. The \texttt{BLDP1} formulation is rarely the best performing formulation, while \texttt{BLDP1+} performs surprisingly well on the most difficult ``Armour`` instances (but rarely otherwise). Overall, the sequence pair base formulation performs the best in a slight majority of the instances (54 of 99). Of this, 35 occur with the sequence pair formulation paired with the new inequalities derived using the embedding approach in this work. The new refined unary formulation is the top performer in roughly a quarter of the benchmark instances (23 of 99), most often augmented with the new inequalities derived herein.

We observe that, although there is no clearly superior formulation or approach across all benchmark instances, it is often the case that one approach will be the clear choice for a particular instance or family of instances. For example, for the \texttt{Bozer97-2-0.0(4)} instance, we observe that none of the existing methods were unable to prove optimality within 7\% relative gap within 4 hours. However, the sequence pair formulation augmented with the newly derived inequalities is able to prove optimality in less than 30 minutes. Additionally, in the \texttt{hp} family of instances we see solve time speed-ups of multiple factors using the new refined unary formulation and valid inequalities over existing formulation methods.

\newcommand{\p}[1]{\tiny (#1\%)}

\begin{table}
    \centering
    {\tiny
\begin{tabular}{r|cccccc|ccccc}
Instance & U & U+ & BLDP1 & BLDP1+ & SP & SP+ & SP+VI & SP+VI3 & RU & RU+VI & RU+VI3\\\hline
\texttt{hp-0.0(4)} & 3067 & 3915 & 7772 & 2873 & 2754 & 2761 & 1463 & 2136 & 1578 & 840 & \cellcolor{blue!25} 656 \\
\texttt{hp-0.0(5)} & 765 & 752 & 928 & 858 & 1363 & 1345 & 836 & 642 & 474 & \cellcolor{blue!25} 350 & 813 \\
\texttt{hp-0.0(6)} & 601 & 357 & 230 & 562 & 427 & 422 & 408 & 384 & 194 & \cellcolor{blue!25} 138 & 207 \\
\texttt{hp-0.1(4)} & 3294 & 3506 & 3941 & 5567 & 4115 & 3977 & 1530 & 2332 & 1471 & \cellcolor{blue!25} 1427 & 1687 \\
\texttt{hp-0.1(5)} & 1076 & 1216 & 619 & 826 & 1230 & 1240 & 548 & 636 & 436 & \cellcolor{blue!25} 304 & 370 \\
\texttt{hp-0.1(6)} & 1198 & 275 & 639 & 349 & 655 & 638 & 346 & 427 & 273 & \cellcolor{blue!25} 157 & 163 \\
\texttt{hp-0.2(4)} & 624 & 1019 & 549 & 1906 & 460 & 437 & 266 & 326 & 279 & 230 & \cellcolor{blue!25} 222 \\
\texttt{hp-0.2(5)} & 808 & 653 & 354 & 1182 & 700 & 708 & 394 & 503 & 470 & \cellcolor{blue!25} 301 & 420 \\
\texttt{hp-0.2(6)} & 504 & 346 & 203 & 345 & 613 & 604 & 302 & 325 & 201 & 350 & \cellcolor{blue!25} 158 \\
\hline
\texttt{apte-0.0(4)} & 4722 & 1040 & \cellcolor{blue!25} 308 & 800 & 414 & 411 & 578 & 962 & 1456 & 1412 & 1433 \\
\texttt{apte-0.0(5)} & 3311 & 571 & 531 & 599 & \cellcolor{blue!25} 165 & 165 & 218 & 198 & 1111 & 398 & 344 \\
\texttt{apte-0.0(6)} & 438 & 352 & 107 & 282 & \cellcolor{blue!25} 90 & 90 & 117 & 169 & 1612 & 340 & 257 \\
\texttt{apte-0.1(4)} & 1935 & 1538 & 810 & 2529 & \cellcolor{blue!25} 247 & 247 & 581 & 343 & 2123 & 1297 & 1057 \\
\texttt{apte-0.1(5)} & 1042 & 250 & \cellcolor{blue!25} 142 & 470 & 201 & 202 & 163 & 285 & 424 & 383 & 215 \\
\texttt{apte-0.1(6)} & 1339 & 478 & 150 & 171 & 205 & 205 & 114 & \cellcolor{blue!25} 86 & 127 & 500 & 158 \\
\texttt{apte-0.2(4)} & 1822 & 883 & 710 & 816 & 694 & 693 & \cellcolor{blue!25} 215 & 237 & 1059 & 528 & 387 \\
\texttt{apte-0.2(5)} & 217 & 184 & 96 & 153 & 102 & 102 & \cellcolor{blue!25} 96 & 103 & 198 & 364 & 159 \\
\texttt{apte-0.2(6)} & 203 & 174 & 62 & 303 & \cellcolor{blue!25} 42 & 42 & 75 & 88 & 149 & 129 & 91 \\
\hline
\texttt{xerox-0.0(4)} & 11467 & 4232 & \tiny{1.96}\% & \tiny{1.92}\% & 737 & \cellcolor{blue!25} 735 & 1022 & 2099 & 3089 & 2639 & 3312 \\
\texttt{xerox-0.0(5)} & 1830 & 1212 & 907 & 896 & 385 & 383 & 500 & 292 & \cellcolor{blue!25} 252 & 1149 & 641 \\
\texttt{xerox-0.0(6)} & 503 & 381 & 2588 & 576 & \cellcolor{blue!25} 167 & 168 & 220 & 324 & 204 & 294 & 331 \\
\texttt{xerox-0.1(4)} & 9264 & 7460 & \tiny{19.99}\% & 10859 & 949 & \cellcolor{blue!25} 948 & 1505 & 1523 & 5642 & 3539 & 7063 \\
\texttt{xerox-0.1(5)} & 778 & 1342 & 1647 & 1016 & 670 & 666 & \cellcolor{blue!25} 254 & 415 & 463 & 685 & 522 \\
\texttt{xerox-0.1(6)} & 481 & 326 & 1047 & 872 & 244 & 243 & \cellcolor{blue!25} 157 & 208 & 323 & 215 & 254 \\
\texttt{xerox-0.2(4)} & 1279 & 1782 & 4331 & 2540 & 1937 & 1934 & \cellcolor{blue!25} 433 & 859 & 793 & 1232 & 2124 \\
\texttt{xerox-0.2(5)} & 1453 & 764 & 1290 & 705 & 783 & 787 & 331 & 420 & \cellcolor{blue!25} 305 & 408 & 435 \\
\texttt{xerox-0.2(6)} & 606 & 261 & 158 & 378 & 123 & 123 & \cellcolor{blue!25} 77 & 94 & 141 & 174 & 192 \\
\hline
\texttt{Camp91-0.0(4)} & 4526 & 7481 & 2150 & 7507 & 931 & \cellcolor{blue!25} 837 & 852 & 1449 & 5498 & 6181 & 3126 \\
\texttt{Camp91-0.0(5)} & 2510 & 4053 & 2553 & 1213 & 429 & \cellcolor{blue!25} 426 & 490 & 499 & 1856 & 2372 & 4513 \\
\texttt{Camp91-0.0(6)} & 1653 & 471 & 793 & 720 & \cellcolor{blue!25} 145 & 148 & 258 & 298 & 721 & 836 & 452 \\
\texttt{Camp91-0.1(4)} & \tiny{3.45}\% & \tiny{4.79}\% & 12967 & \tiny{11.59}\% & 1905 & \cellcolor{blue!25} 1848 & 2064 & 2101 & 6893 & 10669 & \tiny{11.13}\% \\
\texttt{Camp91-0.1(5)} & 9259 & 7702 & \tiny{9.55}\% & \tiny{11.27}\% & 1424 & 1422 & \cellcolor{blue!25} 1144 & 1153 & 3273 & 11078 & \tiny{14.92}\% \\
\texttt{Camp91-0.1(6)} & 660 & 4133 & 5769 & 578 & \cellcolor{blue!25} 181 & 328 & 218 & 396 & 694 & 2271 & 2983 \\
\texttt{Camp91-0.2(4)} & 8584 & 3906 & 2681 & 4109 & \cellcolor{blue!25} 913 & 914 & 1063 & 916 & 3537 & 6301 & 13374 \\
\texttt{Camp91-0.2(5)} & 1259 & 2017 & 1955 & 4134 & \cellcolor{blue!25} 289 & 315 & 310 & 349 & 1795 & 5050 & 4576 \\
\texttt{Camp91-0.2(6)} & 13553 & 11034 & \tiny{12.29}\% & 7393 & 1166 & \cellcolor{blue!25} 1125 & 2259 & 1199 & 4516 & 8845 & 12353 \\
\hline
\texttt{Bozer97-1-0.0(4)} & 1826 & 2046 & 729 & 1512 & 607 & 609 & \cellcolor{blue!25} 354 & 1077 & 1471 & 1456 & 2158 \\
\texttt{Bozer97-1-0.0(5)} & 1159 & 4287 & 1467 & 1059 & 433 & 435 & \cellcolor{blue!25} 411 & 698 & 607 & 1240 & 2306 \\
\texttt{Bozer97-1-0.0(6)} & 1112 & 2050 & 6804 & 2178 & \cellcolor{blue!25} 420 & 420 & 588 & 734 & 494 & 2826 & 6938 \\
\texttt{Bozer97-1-0.1(4)} & 3380 & 2432 & 2482 & 1983 & 896 & \cellcolor{blue!25} 896 & 1059 & 1430 & 6056 & 4190 & 3410 \\
\texttt{Bozer97-1-0.1(5)} & 999 & 1337 & 1075 & 1098 & 312 & \cellcolor{blue!25} 311 & 667 & 1369 & 1082 & 590 & 3450 \\
\texttt{Bozer97-1-0.1(6)} & 2010 & 2503 & 2762 & 3733 & 861 & \cellcolor{blue!25} 850 & 997 & 889 & 2417 & 3352 & 2820 \\
\texttt{Bozer97-1-0.2(4)} & 1661 & 1052 & 610 & 1218 & \cellcolor{blue!25} 312 & 320 & 374 & 402 & 1105 & 1051 & 1836 \\
\texttt{Bozer97-1-0.2(5)} & 1522 & 4009 & 2115 & 1631 & 1253 & 1158 & \cellcolor{blue!25} 740 & 1027 & 1423 & 3272 & 2885 \\
\texttt{Bozer97-1-0.2(6)} & 1271 & 2639 & 1059 & 1322 & \cellcolor{blue!25} 308 & 316 & 401 & 352 & 598 & 3218 & 923 \\
\hline
\texttt{Bozer97-2-0.0(4)} & \tiny{7.05}\% & \tiny{13.43}\% & \tiny{9.21}\% & \tiny{10.63}\% & \tiny{18.45}\% & \tiny{18.43}\% & \cellcolor{blue!25} 1721 & 3818 & \tiny{4.43}\% & \tiny{4.28}\% & \tiny{31.93}\% \\
\texttt{Bozer97-2-0.0(5)} & \tiny{18.16}\% & \tiny{10.15}\% & \tiny{29.32}\% & \tiny{7.71}\% & 12538 & 12654 & \cellcolor{blue!25} 6212 & 6859 & \tiny{9.94}\% & \tiny{2.81}\% & \tiny{3.95}\% \\
\texttt{Bozer97-2-0.0(6)} & \tiny{12.22}\% & 12685 & \tiny{40.59}\% & \tiny{10.32}\% & 4953 & 4913 & 3771 & \cellcolor{blue!25} 3024 & \tiny{0.34}\% & 5491 & 6930 \\
\texttt{Bozer97-2-0.1(4)} & \tiny{17.74}\% & \tiny{17.01}\% & \tiny{28.42}\% & \tiny{20.02}\% & \tiny{13.67}\% & \tiny{13.65}\% & \tiny{8.22}\% & \cellcolor{blue!25} \tiny{7.00}\% & \tiny{15.41}\% & \tiny{20.16}\% & \tiny{28.05}\% \\
\texttt{Bozer97-2-0.1(5)} & \tiny{26.60}\% & \tiny{12.05}\% & \tiny{30.28}\% & \tiny{8.80}\% & \tiny{1.11}\% & \tiny{4.00}\% & 11924 & \cellcolor{blue!25} 9165 & \tiny{10.10}\% & \tiny{3.78}\% & \tiny{15.74}\% \\
\texttt{Bozer97-2-0.1(6)} & \tiny{13.61}\% & \tiny{6.62}\% & \tiny{24.74}\% & \tiny{21.74}\% & 11093 & \tiny{7.38}\% & \cellcolor{blue!25} 1890 & 2047 & 11126 & \tiny{9.29}\% & \tiny{11.97}\% \\
\texttt{Bozer97-2-0.2(4)} & \tiny{4.68}\% & \tiny{1.56}\% & \tiny{6.31}\% & \tiny{3.64}\% & 5264 & 7784 & 5767 & \cellcolor{blue!25} 3607 & 10052 & 10878 & \tiny{4.84}\% \\
\texttt{Bozer97-2-0.2(5)} & \tiny{32.92}\% & \tiny{11.56}\% & \tiny{16.78}\% & \tiny{14.07}\% & \tiny{10.80}\% & \tiny{10.40}\% & \tiny{7.45}\% & \cellcolor{blue!25} \tiny{4.20}\% & \tiny{8.84}\% & \tiny{10.07}\% & \tiny{14.35}\% \\
\texttt{Bozer97-2-0.2(6)} & \tiny{30.67}\% & \tiny{28.28}\% & \tiny{39.68}\% & \tiny{43.09}\% & \tiny{14.36}\% & \tiny{18.19}\% & \tiny{11.62}\% & \cellcolor{blue!25} \tiny{5.90}\% & \tiny{22.58}\% & \tiny{20.05}\% & \tiny{28.53}\% \\
\end{tabular}
}
    \caption{Solution time (or relative gap after 4 hours, if not solved to optimality) for the first 6 benchmark instance families. The first grouping contains approaches from the literature; the second contains approaches using some component (formulation or inequalities) from this work. The best approach for each benchmark is highlighted in blue.}
    \label{table:solve-time-1}
\end{table}

\begin{table}
    \centering
    {\tiny
\begin{tabular}{r|cccccc|ccccc}
Instance & U & U+ & BLDP1 & BLDP1+ & SP & SP+ & SP+VI & SP+VI3 & RU & RU+VI & RU+VI3\\\hline
\texttt{Bazaraa75-1-0.0(4)} & \tiny{1.07}\% & 8017 & 8960 & 7093 & 4229 & 4170 & \cellcolor{blue!25} 3408 & ME & 10025 & 7495 & 10357 \\
\texttt{Bazaraa75-1-0.0(5)} & 12235 & \tiny{5.06}\% & 4257 & ME & 4638 & 4640 & 4000 & 3838 & 4921 & \cellcolor{blue!25} 3211 & 7230 \\
\texttt{Bazaraa75-1-0.0(6)} & 8946 & 9398 & 7645 & \tiny{1.05}\% & 2313 & 2134 & 2010 & 7021 & \cellcolor{blue!25} 1458 & 4583 & 2788 \\
\texttt{Bazaraa75-1-0.1(4)} & 10928 & 8468 & 7371 & \tiny{3.83}\% & 10557 & 10657 & \cellcolor{blue!25} 3938 & 4471 & \tiny{6.85}\% & 13673 & 10898 \\
\texttt{Bazaraa75-1-0.1(5)} & \tiny{5.24}\% & \tiny{3.68}\% & \tiny{1.31}\% & ME & 10393 & \tiny{1.75}\% & \cellcolor{blue!25} 8476 & 11212 & \tiny{3.34}\% & \tiny{4.83}\% & \tiny{2.39}\% \\
\texttt{Bazaraa75-1-0.1(6)} & \tiny{3.49}\% & \tiny{3.10}\% & \tiny{8.64}\% & \tiny{8.07}\% & ME & ME & 13921 & \tiny{100.00}\% & \cellcolor{blue!25} 9195 & \tiny{0.01}\% & \tiny{1.61}\% \\
\texttt{Bazaraa75-1-0.2(4)} & \tiny{5.49}\% & \tiny{2.87}\% & \tiny{4.60}\% & \tiny{9.57}\% & \cellcolor{blue!25} 10215 & 10386 & \tiny{4.44}\% & ME & 11508 & \tiny{8.50}\% & \tiny{11.72}\% \\
\texttt{Bazaraa75-1-0.2(5)} & \tiny{5.09}\% & \tiny{21.29}\% & \tiny{4.05}\% & \tiny{6.58}\% & 5903 & \cellcolor{blue!25} 4894 & 10974 & 10662 & \tiny{1.23}\% & \tiny{3.97}\% & 13139 \\
\texttt{Bazaraa75-1-0.2(6)} & 2516 & 3979 & 1970 & 7441 & 1481 & 1595 & 1262 & 1370 & \cellcolor{blue!25} 1179 & 1626 & 1788 \\
\hline
\texttt{Bazaraa75-2-0.0(4)} & \tiny{37.27}\% & \tiny{38.21}\% & \tiny{26.15}\% & \tiny{30.56}\% & ME & ME & \tiny{25.87}\% & \cellcolor{blue!25} \tiny{21.97}\% & \tiny{51.53}\% & \tiny{54.16}\% & \tiny{37.23}\% \\
\texttt{Bazaraa75-2-0.0(5)} & \tiny{44.73}\% & \tiny{32.44}\% & \cellcolor{blue!25} \tiny{25.82}\% & \tiny{41.09}\% & \tiny{27.48}\% & \tiny{28.74}\% & \tiny{31.71}\% & ME & \tiny{33.85}\% & \tiny{31.16}\% & \tiny{34.83}\% \\
\texttt{Bazaraa75-2-0.0(6)} & \tiny{36.75}\% & \tiny{32.27}\% & \cellcolor{blue!25} \tiny{22.57}\% & \tiny{37.49}\% & ME & ME & \tiny{29.15}\% & \tiny{35.82}\% & \tiny{24.80}\% & ME & \tiny{37.96}\% \\
\texttt{Bazaraa75-2-0.1(4)} & \tiny{42.35}\% & \tiny{41.74}\% & \tiny{31.43}\% & \tiny{42.17}\% & ME & ME & ME & ME & \tiny{35.42}\% & \tiny{31.20}\% & \cellcolor{blue!25} \tiny{29.19}\% \\
\texttt{Bazaraa75-2-0.1(5)} & \tiny{33.37}\% & \tiny{40.72}\% & \tiny{44.03}\% & \tiny{28.60}\% & ME & ME & \cellcolor{blue!25} \tiny{25.77}\% & ME & ME & ME & \tiny{33.84}\% \\
\texttt{Bazaraa75-2-0.1(6)} & \tiny{30.08}\% & \tiny{37.95}\% & \cellcolor{blue!25} \tiny{20.02}\% & \tiny{56.93}\% & \tiny{31.67}\% & \tiny{31.66}\% & \tiny{29.06}\% & ME & \tiny{24.24}\% & \tiny{25.15}\% & \tiny{24.35}\% \\
\texttt{Bazaraa75-2-0.2(4)} & \tiny{42.53}\% & \tiny{27.64}\% & \tiny{27.79}\% & \cellcolor{blue!25} \tiny{26.52}\% & ME & ME & ME & ME & \tiny{35.63}\% & \tiny{35.75}\% & \tiny{31.02}\% \\
\texttt{Bazaraa75-2-0.2(5)} & \tiny{37.17}\% & \tiny{40.36}\% & \cellcolor{blue!25} \tiny{29.01}\% & \tiny{37.53}\% & \tiny{29.88}\% & \tiny{29.42}\% & ME & \tiny{29.81}\% & \tiny{29.53}\% & \tiny{29.20}\% & \tiny{32.36}\% \\
\texttt{Bazaraa75-2-0.2(6)} & \tiny{29.95}\% & \tiny{37.58}\% & \tiny{32.22}\% & \tiny{34.36}\% & \tiny{24.93}\% & ME & ME & ME & \cellcolor{blue!25} \tiny{20.48}\% & ME & \tiny{32.60}\% \\
\hline
\texttt{Bozer91-0.0(4)} & \tiny{45.60}\% & \tiny{37.60}\% & \tiny{63.39}\% & \tiny{32.56}\% & \tiny{39.41}\% & \tiny{39.41}\% & \cellcolor{blue!25} \tiny{20.75}\% & \tiny{24.63}\% & \tiny{45.95}\% & \tiny{38.62}\% & \tiny{38.30}\% \\
\texttt{Bozer91-0.0(5)} & \tiny{55.81}\% & \tiny{30.63}\% & \tiny{60.98}\% & \tiny{22.83}\% & \tiny{32.12}\% & \tiny{31.99}\% & \tiny{25.10}\% & \cellcolor{blue!25} \tiny{21.99}\% & \tiny{41.39}\% & \tiny{37.89}\% & \tiny{29.61}\% \\
\texttt{Bozer91-0.0(6)} & \tiny{100.00}\% & \tiny{42.60}\% & \tiny{51.66}\% & \cellcolor{blue!25} \tiny{23.86}\% & \tiny{35.56}\% & \tiny{34.90}\% & \tiny{24.71}\% & \tiny{25.52}\% & \tiny{25.26}\% & \tiny{31.95}\% & \tiny{32.91}\% \\
\texttt{Bozer91-0.1(4)} & \tiny{43.60}\% & \tiny{40.54}\% & \tiny{52.11}\% & \tiny{36.46}\% & \tiny{36.61}\% & \tiny{36.62}\% & \cellcolor{blue!25} \tiny{24.83}\% & \tiny{24.87}\% & \tiny{42.88}\% & \tiny{35.57}\% & \tiny{32.48}\% \\
\texttt{Bozer91-0.1(5)} & \tiny{41.24}\% & \tiny{43.68}\% & \tiny{52.88}\% & \tiny{31.66}\% & \tiny{33.32}\% & \tiny{31.13}\% & \cellcolor{blue!25} \tiny{22.80}\% & \tiny{22.85}\% & \tiny{46.12}\% & \tiny{34.17}\% & \tiny{30.89}\% \\
\texttt{Bozer91-0.1(6)} & \tiny{47.96}\% & \tiny{33.48}\% & \tiny{60.68}\% & \cellcolor{blue!25} \tiny{21.43}\% & \tiny{34.83}\% & \tiny{36.51}\% & \tiny{22.30}\% & \tiny{25.02}\% & \tiny{34.69}\% & \tiny{31.16}\% & \tiny{30.98}\% \\
\texttt{Bozer91-0.2(4)} & \tiny{49.25}\% & \tiny{38.45}\% & \tiny{73.86}\% & \tiny{29.83}\% & \tiny{38.49}\% & \tiny{38.42}\% & \tiny{22.09}\% & \cellcolor{blue!25} \tiny{21.36}\% & \tiny{59.05}\% & \tiny{29.32}\% & \tiny{35.00}\% \\
\texttt{Bozer91-0.2(5)} & \tiny{41.32}\% & \tiny{30.72}\% & \tiny{46.12}\% & \tiny{26.24}\% & \tiny{28.47}\% & \tiny{27.75}\% & \tiny{16.16}\% & \cellcolor{blue!25} \tiny{15.00}\% & \tiny{34.63}\% & \tiny{35.58}\% & \tiny{30.08}\% \\
\texttt{Bozer91-0.2(6)} & \tiny{48.27}\% & \tiny{34.59}\% & \tiny{53.13}\% & \tiny{41.43}\% & \tiny{28.21}\% & \tiny{31.67}\% & \tiny{26.72}\% & \cellcolor{blue!25} \tiny{14.24}\% & \tiny{39.77}\% & \tiny{34.65}\% & \tiny{32.50}\% \\
\hline
\texttt{Armour62-1-0.0(4)} & \tiny{100.00}\% & \tiny{62.31}\% & \tiny{72.30}\% & \cellcolor{blue!25} \tiny{60.25}\% & ME & ME & ME & ME & \tiny{68.95}\% & \tiny{62.03}\% & \tiny{62.24}\% \\
\texttt{Armour62-1-0.0(5)} & \tiny{100.00}\% & \tiny{100.00}\% & \tiny{77.77}\% & \cellcolor{blue!25} \tiny{64.19}\% & ME & ME & ME & ME & \tiny{69.05}\% & \tiny{69.60}\% & \tiny{65.03}\% \\
\texttt{Armour62-1-0.0(6)} & \tiny{100.00}\% & \tiny{100.00}\% & \tiny{70.20}\% & \cellcolor{blue!25} \tiny{66.43}\% & ME & ME & ME & ME & \tiny{67.50}\% & \tiny{70.41}\% & ME \\
\texttt{Armour62-1-0.1(4)} & \tiny{73.88}\% & \tiny{100.00}\% & \tiny{69.54}\% & \cellcolor{blue!25} \tiny{59.93}\% & ME & ME & ME & ME & \tiny{71.93}\% & \tiny{64.70}\% & \tiny{62.85}\% \\
\texttt{Armour62-1-0.1(5)} & \tiny{100.00}\% & \tiny{100.00}\% & \tiny{76.18}\% & \tiny{66.73}\% & ME & ME & ME & ME & \tiny{75.34}\% & \tiny{68.35}\% & \cellcolor{blue!25} \tiny{60.94}\% \\
\texttt{Armour62-1-0.1(6)} & \tiny{100.00}\% & \tiny{100.00}\% & \tiny{75.32}\% & \cellcolor{blue!25} \tiny{66.38}\% & ME & ME & ME & ME & \tiny{67.81}\% & \tiny{68.87}\% & \tiny{67.24}\% \\
\texttt{Armour62-1-0.2(4)} & \tiny{100.00}\% & \tiny{100.00}\% & \tiny{71.30}\% & \cellcolor{blue!25} \tiny{60.06}\% & ME & ME & ME & ME & \tiny{71.66}\% & \tiny{62.84}\% & \tiny{63.68}\% \\
\texttt{Armour62-1-0.2(5)} & \tiny{100.00}\% & \tiny{100.00}\% & \tiny{83.28}\% & \tiny{64.84}\% & ME & ME & ME & ME & \tiny{73.09}\% & \tiny{65.92}\% & \cellcolor{blue!25} \tiny{59.36}\% \\
\texttt{Armour62-1-0.2(6)} & \tiny{100.00}\% & \tiny{100.00}\% & \tiny{87.65}\% & \cellcolor{blue!25} \tiny{59.76}\% & ME & ME & ME & ME & \tiny{75.79}\% & \tiny{73.12}\% & \tiny{65.27}\% \\
\hline
\texttt{Armour62-2-0.0(4)} & \tiny{100.00}\% & \tiny{60.68}\% & \tiny{63.44}\% & \cellcolor{blue!25} \tiny{60.04}\% & ME & ME & ME & ME & \tiny{64.23}\% & \tiny{62.93}\% & \tiny{64.38}\% \\
\texttt{Armour62-2-0.0(5)} & \tiny{100.00}\% & \tiny{100.00}\% & \tiny{69.42}\% & \cellcolor{blue!25} \tiny{60.41}\% & ME & ME & ME & \tiny{74.11}\% & \tiny{70.65}\% & \tiny{63.24}\% & \tiny{66.76}\% \\
\texttt{Armour62-2-0.0(6)} & \tiny{100.00}\% & \tiny{100.00}\% & \tiny{72.00}\% & \cellcolor{blue!25} \tiny{58.22}\% & ME & ME & \tiny{70.84}\% & \tiny{63.48}\% & \tiny{71.49}\% & \tiny{70.60}\% & \tiny{61.51}\% \\
\texttt{Armour62-2-0.1(4)} & \tiny{70.13}\% & \tiny{100.00}\% & \tiny{67.45}\% & \tiny{63.29}\% & ME & ME & ME & ME & \tiny{67.06}\% & \tiny{65.71}\% & \cellcolor{blue!25} \tiny{60.74}\% \\
\texttt{Armour62-2-0.1(5)} & \tiny{100.00}\% & \tiny{100.00}\% & \tiny{67.01}\% & \tiny{61.39}\% & ME & ME & ME & ME & \tiny{71.38}\% & \tiny{68.30}\% & \cellcolor{blue!25} \tiny{60.20}\% \\
\texttt{Armour62-2-0.1(6)} & \tiny{100.00}\% & \tiny{100.00}\% & \tiny{81.03}\% & \cellcolor{blue!25} \tiny{55.84}\% & ME & ME & \tiny{73.37}\% & ME & \tiny{74.06}\% & \tiny{68.71}\% & \tiny{62.16}\% \\
\texttt{Armour62-2-0.2(4)} & \tiny{100.00}\% & \tiny{100.00}\% & \tiny{75.49}\% & \cellcolor{blue!25} \tiny{62.18}\% & ME & ME & ME & ME & \tiny{70.81}\% & \tiny{64.45}\% & ME \\
\texttt{Armour62-2-0.2(5)} & \tiny{100.00}\% & \tiny{100.00}\% & \tiny{75.65}\% & \tiny{66.07}\% & ME & ME & ME & ME & \tiny{70.76}\% & \tiny{67.01}\% & \cellcolor{blue!25} \tiny{62.35}\% \\
\texttt{Armour62-2-0.2(6)} & \tiny{100.00}\% & \tiny{100.00}\% & \tiny{90.82}\% & \cellcolor{blue!25} \tiny{60.32}\% & ME & ME & ME & ME & \tiny{74.19}\% & \tiny{69.65}\% & \tiny{64.48}\% \\
\end{tabular}
}
    \caption{Solution time (or relative gap after 4 hours, if not solved to optimality) for the last 5 benchmark instance families. The first grouping contains approaches from the literature; the second contains approaches using some component (formulation or inequalities) from this work. The best approach for each benchmark is highlighted in blue.}
    \label{table:solve-time-2}
\end{table}

The results in Tables \ref{table:solve-time-1} and \ref{table:solve-time-2} suggests that none of the approaches herein will be a clear winner on all instances, but that some combination of them can be used to tackle difficult problems. In particular, we recommend trying both the sequence pair and refined unary formulations, along with some subsets of the inequalities derived in this work to solve the FLP instance most efficiently.

\section{Conclusion}\label{sec:conclusions}
In this work, we  presented a case study on systematically building strong formulations for disjunctive sets; namely, for the floor layout problem. We  used the embedding approach of \cite{Vielma:2015a} to generate MIP formulations of the FLP and have observed how, by varying our inputs to the procedure, we are able to reconstruct all existing MIP formulations for the problem, produce new formulations, and discover valid inequalities. We also showed how valid inequalities generated for one formulation can often be translated to seemingly unrelated formulations. Finally, we presented computational results showing how the developed techniques can been used to solve previously unsolved benchmark instances.

\section*{Acknowledgments}
    This material is based upon work supported by the National Science Foundation Graduate Research Fellowship under Grant No. 1122374 and Grant CMMI-1351619.

\bibliographystyle{ormsv080}
{\footnotesize \bibliography{master}}

\newpage
\appendix

\section{Standard extended formulation for FLP}\label{extendedformulationsec}
The following corollary gives an extended formulation for $\hat{\scrL}_{i,j}$ using a standard approach by Balas, Jeroslow and Lowe.
\begin{corollary}\label{balascoro}
	The following is an ideal extended formulation for $\hat{\scrL}_{i,j}$:
	\begin{subequations}\label{eqn:balasform}
	\begin{alignat}{3}
		\frac{1}{2}\ell^s_{p,q} \leq c^s_{p,q} &\leq L^sv_q - \frac{1}{2}\ell^s_{p,q} &\quad& \forall s \in \{x,y\}, p \in \{i,j\}, q \in \llbracket 4 \rrbracket \\
		lb^s_pv_q \leq \ell^s_{p,q} &\leq ub^s_pv_q &\quad& \forall p \in \{i,j\}, q \in \llbracket 4 \rrbracket \\
		c^y_{i,1} - c^y_{j,1} + \frac{1}{2}(\ell^y_{i,1} + \ell^y_{j,1}) &\leq 0 \\
		c^x_{i,2} - c^x_{j,2} + \frac{1}{2}(\ell^x_{i,2} + \ell^x_{j,2}) &\leq 0 \\
		c^y_{j,3} - c^y_{i,3} + \frac{1}{2}(\ell^y_{i,3} + \ell^y_{j,3}) &\leq 0 \\
		c^x_{j,4} - c^x_{i,4} + \frac{1}{2}(\ell^x_{i,4} + \ell^x_{j,4}) &\leq 0 \\
		\sum_{i=1}^4 c^s_{p,i} &= c^s_p &\quad& \forall s \in \{x,y\}, p \in \{i,j\} \\
		\sum_{i=1}^4 \ell^s_{p,i} &= \ell^s_p &\quad& \forall s \in \{x,y\}, p \in \{i,j\} \\
		\sum_{i=1}^4 v_i &= 1 \label{int1}\\
		v &\in \{0,1\}^4.\label{int2}
	\end{alignat}
	\end{subequations}
\end{corollary}
\proof{}
	This follows from Proposition 4.2 in  \cite{Vielma:2015}.
\endproof

\section{Proof for Theorem \ref{thm:unary}} \label{app:unary-proof}
\proof{}
    First, we note that \eqref{eqn:unary-formulation} is just formulation \eqref{eqn:unary-bigM} with a tightened form of constraints \eqref{eqn:tight-sitb}.
	For validity, first we want to show that the stay-on-the-floor constraints can be tightened to \eqref{eqn:tight-sitb} by the following case analysis. Consider $s=x$, $p=i$, and $q = j$; the other constraints follow analogously.
	\begin{itemize}
		\item \underline{$u^x_{i,j}=0, u^x_{j,i}=0$} Reduces to the linear constraints in \eqref{eqn:sitb} defining $Q^{lb}$.
		\item \underline{$u^x_{i,j}=1, u^x_{j,i}=0$} The first inequality is unchanged. For the second, we note that $u^x_{i,j} = 1$ means that
			\begin{align*}
				c^x_i + \frac{1}{2}\ell^x_i &\leq c^y_j - \frac{1}{2}\ell^x_j \\
				&\leq \left(L^x-\frac{1}{2}\ell^x_j\right) - \frac{1}{2}\ell^x_j \quad \text{from \eqref{eqn:sitb}} \\
				&= L^x - \ell^x_j \\
				&\leq L^x - lb^x_j
			\end{align*}
		\item \underline{$u^x_{i,j}=0, u^x_{j,i}=1$} Same argument as the case before, but with the branch $c^x_j + \frac{1}{2}\ell^x_j \leq c^y_i - \frac{1}{2}\ell^x_i$ and the first inequality in \eqref{eqn:tight-sitb}.
		\item \underline{$u^x_{i,j}=1, u^x_{j,i}=1$} Not feasible by \eqref{eqn:sum-to-one}.
	\end{itemize}

	To show the formulation is ideal, we want to show that every extreme point of its relaxation has integral value for $u$.

	Consider the following system, which can be thought of as the projection of the relaxation of \eqref{eqn:unary-formulation} onto just the $y$ variables (the argument for $x$ variables is identical). Under the assumption that $lb^y_i + lb^y_j < L^y$, the following is full-dimensional:
	\begin{subequations} \label{eqn:one-dimensional-unary}
	\begin{gather}
		\frac{1}{2}\ell_i + lb^y_j u_{j,i} \leq c_i \label{eqn:42} \\
		c_i \leq L - \frac{1}{2}\ell_i - lb^y_ju_{i,j} \label{eqn:43} \\
		\frac{1}{2}\ell_j + lb^y_i u_{i,j} \leq c_j \label{eqn:44} \\
		c_j \leq L - \frac{1}{2}\ell_j - lb^y_iu_{j,i} \label{eqn:45} \\
		c_i + \frac{1}{2}\ell_i \leq c_j - \frac{1}{2}\ell_j + L^y(1-u_{i,j}) \label{eqn:46} \\
		c_j + \frac{1}{2}\ell_j \leq c_i - \frac{1}{2}\ell_i + L^y(1-u_{j,i}) \label{eqn:47} \\
		\ell_i \geq lb^y_i \label{eqn:48} \\
		\ell_j \geq lb^y_j \label{eqn:49} \\
		u_{i,j} \geq 0 \label{eqn:50} \\
		u_{i,j} \leq 1 \label{eqn:51} \\
		u_{j,i} \geq 0 \label{eqn:52} \\
		u_{j,i} \leq 1 \label{eqn:53} \\
		u_{i,j} + u_{j,i} \leq 1 \label{eqn:54}.
	\end{gather}
	\end{subequations}

	Take some arbitrary feasible $(\hat{c},\hat{\ell},\hat{u})$ where $0 < \hat{u}_{i,j} < 1$ is fractional. We wish to show that this is not an extreme point. The argument for fractional $\hat{u}_{j,i}$ follows in the same way. We will consider a partition of all possible cases in the following way:

	\begin{enumerate}
		\item \eqref{eqn:52} is active
			\begin{enumerate}
				\item \eqref{eqn:47} is active
				\item \eqref{eqn:47} is not active
			\end{enumerate}
		\item \eqref{eqn:52} is not active
			\begin{enumerate}
				\item \eqref{eqn:46} and \eqref{eqn:47} are both active
					\begin{enumerate}
						\item \eqref{eqn:54} is active
						\item \eqref{eqn:54} is not active
					\end{enumerate}
				\item At most one of \eqref{eqn:46} and \eqref{eqn:47} are active
					\begin{enumerate}
						\item \eqref{eqn:54} is not active
						\item \eqref{eqn:54} is active
					\end{enumerate}
			\end{enumerate}
	\end{enumerate}
	Note in particular that, since $\hat{u}_{i,j} > 0$, constraint \eqref{eqn:54} immediately implies that $\hat{u}_{j,i} < 1$ and that (\ref{eqn:50}-\ref{eqn:51}) cannot be active. In each of these cases, we will argue that the solution is not extreme, either because there are not the requisite six constraints active, or because the selection leads to a contradiction.

	Also, note that \eqref{eqn:43} and \eqref{eqn:44} both being active implies
	\[
		\hat{c}_i + \frac{1}{2}\hat{\ell}_i = \hat{c}_j - \frac{1}{2}\hat{\ell}_j + L - (lb_i+lb_j)\hat{u}_{i,j} > \hat{c}_j - \frac{1}{2}\hat{\ell}_j + L(1-\hat{u}_{i,j})
	\]
	under the assumptions on lower bounds ($L > lb_i+lb_j$) and on $\hat{u}_{i,j}$ fractional; similarly for \eqref{eqn:42} and \eqref{eqn:45} together for $\hat{u}_{j,i} > 0$. Therefore, these pairs cannot be active together when these conditions are met. We will use this observation in the following.

	\paragraph{\underline{1.a}}
	First consider the case that \eqref{eqn:47} and \eqref{eqn:52} are active. Clearly \eqref{eqn:53} and \eqref{eqn:54} cannot be active. Then we must have \eqref{eqn:42} and \eqref{eqn:45} active as well: take their sum, which must hold with equality as it is equivalent to \eqref{eqn:47} when $\hat{u}_{j,i}=0$. However, this also implies that at most three of these four active constraints are linearly independent.

	We cannot have \eqref{eqn:43} or \eqref{eqn:44} active under our assumptions on the lower bounds (take their sum with \eqref{eqn:42} or \eqref{eqn:45}, respectively). The only remaining possibility is if \eqref{eqn:46}, \eqref{eqn:48}, and \eqref{eqn:49} are all active. However, then summing \eqref{eqn:46} and \eqref{eqn:47} and reducing leads to
	\[
		lb_i + lb_j = L^y(2-\hat{u}_{i,j}) > L^y,
	\]
	a contradiction. Therefore, the point is not extreme.

	\paragraph{\underline{1.b}}
	Now assume that \eqref{eqn:47} is not active and \eqref{eqn:52} is. This implies that at most one of \eqref{eqn:42} and \eqref{eqn:45} can be active, as their sum is equal to \eqref{eqn:47} (when $\hat{u}_{j,i}=0$) as mentioned in the previous case. From our note above, at most one of \eqref{eqn:43} and \eqref{eqn:44} can be active at once.

	Given this, the only possibility that the point is extreme is if one of \eqref{eqn:42} or \eqref{eqn:45}, one of \eqref{eqn:43} or \eqref{eqn:44}, and \eqref{eqn:46}, \eqref{eqn:48}, and \eqref{eqn:49} are all active. If \eqref{eqn:42} and \eqref{eqn:43} are both active in this setting, their sum implies
	\[
		0 = L - lb_i - lb_j\hat{u}_{i,j} > L - lb_i - lb_j,
	\]
	a contradiction. If \eqref{eqn:42} and \eqref{eqn:44} are simultaneously active in this setting, then the sum $\eqref{eqn:42} - \eqref{eqn:44} + \eqref{eqn:46}$ yields
	\[
		(L - lb_i)(1-\hat{u}_{i,j}) = 0,
	\]
	a contradiction. Therefore, \eqref{eqn:42} is not active. If \eqref{eqn:43} and \eqref{eqn:45} are both active in this setting, then we get a similar result for $\eqref{eqn:43} - \eqref{eqn:45} - \eqref{eqn:46}$. This implies that we can have at most one of (\ref{eqn:42}-\ref{eqn:45}) active, leaving us with too few active constraints for the point to be extreme.

	\paragraph{\underline{2.a.i}}
	Now consider the case that $0 < \hat{u}_{j,i} \leq 1 - \hat{u}_{i,j} < 1$; clearly \eqref{eqn:53} cannot be active. If \eqref{eqn:46}, \eqref{eqn:47}, and \eqref{eqn:54} all are active, they imply
	\begin{equation} \label{eqn:useful-eqn}
		\hat{\ell}_i + \hat{\ell}_j = L^y,
	\end{equation}
	which means that \eqref{eqn:48} and \eqref{eqn:49} cannot both be active by the assumption on lower bounds. If both \eqref{eqn:44} and \eqref{eqn:45} are active together, this implies that \eqref{eqn:48} is active; therefore, \eqref{eqn:49} is not active. Moreover, the discussed active constraints have rank five at this point. By the notes above, \eqref{eqn:44} and \eqref{eqn:45} being active imply that \eqref{eqn:42} and \eqref{eqn:43} cannot be active. The (six) active constraints are linearly dependent, and so the point is not extreme.

	The same argument holds if both \eqref{eqn:42} and \eqref{eqn:43} are active. If \eqref{eqn:42} and \eqref{eqn:44} both are active, then necessarily \eqref{eqn:43} is inactive; summing \eqref{eqn:42} and \eqref{eqn:43} yields
	\[
		\hat{\ell}_i + lb_j < L^y,
	\]
	which says that \eqref{eqn:49} cannot be active from \eqref{eqn:useful-eqn}. Presume then that \eqref{eqn:48} is active, leaving us with \eqref{eqn:42}, \eqref{eqn:44}, \eqref{eqn:46}, \eqref{eqn:47}, \eqref{eqn:48}, and \eqref{eqn:54} active. Taking the resulting description for the extreme point represented by this system of active constraints yields $\hat{u}_{i,j} = \frac{-lb_i-lb_j}{L-lb_i-lb_j} < 0$, a contradiction.

	The same argument holds if \eqref{eqn:43} and \eqref{eqn:45} are both active.

	\paragraph{\underline{2.a.ii}}
	Now assume that \eqref{eqn:54} is not active and that \eqref{eqn:46} and \eqref{eqn:47} are active. Then at most one of \eqref{eqn:42} and \eqref{eqn:43} can be active, else we imply that
	\[
		\hat{u}_{i,j} + \hat{u}_{j,i} = \frac{L-\hat{\ell}_i}{lb_j} \geq \frac{L-lb_i}{lb_j} \geq \frac{lb_j}{lb_j} = 1,
	\]
	which contradicts \eqref{eqn:54} not being active. The same holds for \eqref{eqn:44} and \eqref{eqn:45}. Also, we cannot have both \eqref{eqn:48} and \eqref{eqn:49} active, since along with \eqref{eqn:46} and \eqref{eqn:47} active they imply
	\[
		lb^y_i + lb^y_j = L(2-\hat{u}_{i,j}-\hat{u}_{j,i}) > L.
	\]
	Therefore, at most five constraints are active, and we are not extreme.

	\paragraph{\underline{2.b.i}}
	Now $0 < \hat{u}_{j,i} < 1$ and w.l.o.g. \eqref{eqn:47} is not active. Our statement at the beginning implies that at most two of (\ref{eqn:42}-\ref{eqn:45}) are active. If \eqref{eqn:54} is not active, this leaves at most five active constraints, and so the point is not extreme.

	\paragraph{\underline{2.b.ii}}
	Now assume that \eqref{eqn:47} is not active but \eqref{eqn:54} is. By the argument in 2.b.i, there are at most six constraints (\eqref{eqn:46}, \eqref{eqn:48}, \eqref{eqn:49}, \eqref{eqn:54}, and two of (\ref{eqn:42}-\ref{eqn:45})) that could be active. In particular, \eqref{eqn:46}, \eqref{eqn:48}, and \eqref{eqn:49} must be active at an extreme point, so presume they are.

	Assume for the first case that both \eqref{eqn:42} and \eqref{eqn:44} are active. Computing $\eqref{eqn:42} - \eqref{eqn:44} + \eqref{eqn:46}$ and reducing using the other active constraints yields
	\[
		(lb_i + lb_j - L)\hat{u}_{j,i} = 0,
	\]
	a contradiction. Alternatively, summing \eqref{eqn:42}, \eqref{eqn:45}, and \eqref{eqn:46} yields
	\[
		(lb_i + lb_j - L)(1+\hat{u}_{j,i}) = 0,
	\]
	also a contradiction. Similarly for \eqref{eqn:43}, \eqref{eqn:44}, and \eqref{eqn:46}, as for \eqref{eqn:43}, \eqref{eqn:45}, and \eqref{eqn:46}. This exhausts all possible combinations of active constraints, and so the point is not extreme.

	\paragraph{Piecing together direction-wise formulations}

	Now that we have established that the formulation, when restricted to a single direction, is ideal, it remains to show that the Cartesian product of formulations for both directions, along with the restriction $\{u^x_{i,j}+u^x_{j,i}+u^y_{i,j}+u^y_{j,i}=1\}$, is also ideal. To see this, consider a potential fractional extreme point $(\hat{c},\hat{\ell},\hat{u})$ for the relaxation of the original formulation \eqref{eqn:unary-formulation}. Then at least 12 active constraints at such an extreme point, one of which will be \eqref{eqn:sum-to-one}. Of the 11 remaining that must exist, there will be one direction for which there are at least six active constraints. Consider three cases.

	\begin{enumerate}
		\item The direction with six active constraints has a fractional component in $u$. Then this implies an extreme point for the auxiliary system \eqref{eqn:one-dimensional-unary} with fractional component, a contradiction.
		\item The direction with fractional component (w.l.o.g. $y$) has five active constraints. Then those five active constraints, along with \eqref{eqn:54}, induce a fractional extreme point for the auxiliary system \eqref{eqn:one-dimensional-unary} for $y$, a contradiction.
		\item The direction with fractional component (w.l.o.g. $y$) has fewer than five active constraints. This implies that there are at least seven linearly independent active constraints for the auxiliary system for $x$, a contradiction (since its dimensionality is only six).
	\end{enumerate}
	Therefore, any fractional extreme point induces a fractional extreme point on \eqref{eqn:one-dimensional-unary}, which we have shown is impossible, and so we are done.

\endproof

\section{Big-$M$ formulation for arbitrary encodings and proof of Proposition~\ref{graybigm} and \ref{prop:refined}}

While validity of most formulations in this work can be checked directly, the following generic big-$M$ formulations approach was instrumental for their construction. We include a proof of its validity for completeness.

\begin{theorem} \label{thm:bigM}
	Take $Q$ as a compact convex set, $D = \bigvee_{k=1}^K [A^k x \leq b^k]$, and $K$ distinct vectors $\{v^k\}_{k=1}^K \subseteq \{0,1\}^m$. Take any MIP formulation $V$ for the set $C=\{v^k\}_{k=1}^K$, and some affine functions $R^k_l$ such that
	\[
		R^k_l(v^s) \begin{cases} = b^k_l & k = s \\
 								 \geq \max_{x \in Q(s)} (A^k)_l x & \text{o.w.}
 					 \end{cases} \quad \forall k,l,
	\]
    where $Q(s) \defeq \{x \in Q: A^sx \leq b^s\}$. Then
\begin{subequations}
\begin{gather}
	(x,v) \in Q \times V \\
	(A^k)_l x \leq R^k_l(v) \quad \forall k, l \label{eqn:bigM}
\end{gather}
\end{subequations}
is a valid formulation for $\En(Q,D,C)$ (and, hence, for $\{x \in Q: D\}$).
\end{theorem}
\proof{}
	It is clear from the definition of $D$ that for any $x \in \{x \in Q: D\}$, there is some branch $k$ such that $A^k x \leq b^k$, and so $(x,v^k)$ is feasible for the MIP formulation by the construction of the $R^i_j$. To show that any feasible solution for the MIP formulation lies in $\{x \in Q: D\}$, consider some feasible $(x,v)$ for the MIP formulation.
	Then $x \in Q$ and $v \in V$ implies that $v = v^k$ for some $k$. Then
	\[
		(A^k)_jx \leq b^\ell_j \quad \forall j,
	\]
	implying that $x$ satisfies the corresponding branch $k$ of the disjunction $D$.
\endproof

\subsection{Proof of Proposition \ref{graybigm}} \label{app:gray-bigM}
\proof{}
	Apply Theorem \ref{thm:bigM} with $D^4_{i,j}$, $V = \{0,1\}^2$, and
	\begin{align*}
	    R^1(w) &= L^y(  w_1+w_2) \\
	    R^2(w) &= L^x(1-w_1+w_2) \\
	    R^3(w) &= L^y(2-w_1-w_2) \\
	    R^4(w) &= L^x(1+w_1-w_2).
	\end{align*}
\endproof

\subsection{Proof of Proposition \ref{prop:refined}} \label{app:refined}
\proof{}
	The general proof technique is as follows. First, we will construct the components needed to apply Theorem \ref{thm:bigM}: namely, a ground set describing shared constraints across all feasible layouts, a disjunction we are interested in modeling, a valid formulation for the codes, and some big-$M$ functions that encapsulate the logic between the codes and the branches of the disjunction. This will leave us with a valid formulation for our set $\En(Q^{lb},D^8,C^8)$. We will then do ad-hoc tightening of some of the resulting constraints, giving the system described in \eqref{eqn:refined-unary-formulation}.

	First, we choose the ground set $Q^{lb}$ and the disjunction $D^8$. We see that
	\[
		V \defeq \left\{(z^y_{i,j},z^x_{i,j},z^y_{j,i},z^x_{j,i}) \in \{0,1\}^4: z^x_{i,j}+z^x_{j,i}+z^y_{i,j}+z^y_{j,i} \geq 1, \: z^x_{i,j}+z^x_{j,i} \leq 1, \: z^y_{i,j}+z^y_{j,i} \leq 1\right\}
	\]
	is a valid formulation for $C^8$. Choose big-$M$ functions $R$ based on the disjunction $D^8$ in the following way. Take $T_{f,g}$ as the $g$-th clause defining $br^f$ in Section \ref{sec:refined-disjunction} (recall that ${D^8 = \bigvee_{k=1}^8 br^k}$). for example, $T_{3,2} = \scrB_i \not\leftarrow_y \scrB_j$. Then take
	\[
		R_{f,g}(z) = \begin{cases}
			L^s(1-z^s_{p,q})             & T_{f,g} = \scrB_p \leftarrow_s \scrB_q \\
			(L^s-lb^s_i-lb^s_j)z^s_{p,q} & T_{f,g} = \scrB_p \not\leftarrow_s \scrB_q.
		\end{cases}
	\]

	Note that, when $T_{f,g}$ is a statement of the form ``$\scrB_p$ precedes $\scrB_q$ in direction $s$'', we get the same big-$M$ functions as appeared in the unary formulation for the same logic.

	Now apply Theorem \ref{thm:bigM} and recover the valid formulation $\{(c,\ell,z) \in \bbR^{4+4} \times \{0,1\}^4: (\ref{eqn:refined-untight-sitb}-\ref{eqn:refined-untight-new}),(\ref{eqn:ru-nonoverlap}-\ref{eqn:refined-keeper-last})\}$, where
	\begin{subequations}
	\begin{align}
	    \frac{1}{2}\ell^s_k \leq c^s_k \leq L^s - \frac{1}{2}\ell^s_k \quad &\forall s \in \{x,y\}, k \in \{i,j\} \label{eqn:refined-untight-sitb} \\
	    c^s_p + \frac{1}{2}\ell^s_p + (L^s - lb^s_i - lb^s_j)z^s_{p,q} \geq c^s_q - \frac{1}{2}\ell^s_q \quad &\forall s \in \{x,y\}, \{p,q\} = \{i,j\} \label{eqn:refined-untight-new}.
	\end{align}
	\end{subequations}
	Note that, since the branches of the disjunction share constraints (and corresponding big-$M$ functions $R$), many of the resulting constraints will be equivalent and duplicates can be removed.

	We now wish to tighten some of these constraints by lifting them in an ad-hoc manner. By the same argument as in the proof for Theorem \ref{thm:unary}, we may tighten \eqref{eqn:refined-untight-sitb} to \eqref{eqn:tight-sitb}.

	To tighten the new constraints \eqref{eqn:refined-untight-new}, we can do a case analysis. Consider $s = y$, $p = i$, and $q=j$; the others follow analogously.
	\begin{itemize}
		\item \underline{$z^y_{i,j}=0, z^y_{j,i}=0$} Reduces to the linear constraint $c^y_i + \frac{1}{2}\ell^y_i \geq c^y_j - \frac{1}{2}\ell^y_j$, which follows from $z^y_{i,j} = 0$.
		\item \underline{$z^y_{i,j}=1, z^y_{j,i}=0$} We have that in this case
		\begin{align*}
			c^y_i - \frac{1}{2}\ell^y_i - (c^y_j + \frac{1}{2}\ell^y_j) \geq -L^y
		\end{align*}
		and adding $\ell^y_i + \ell^y_j \geq lb^y_i + lb^y_j$ to this gives the desired inequality
		\[
			c^y_i + \frac{1}{2}\ell^y_i - (c^y_j - \frac{1}{2}\ell^y_j) \geq lb^y_i + lb^y_j - L^y.
		\]
		\item \underline{$z^y_{i,j}=0, z^y_{j,i}=1$} In this case, we have that
		\[
			c^y_i - \frac{1}{2}\ell^y_i - (c^y_j + \frac{1}{2}\ell^y_j) \geq 0.
		\]
		We can add $\ell^y_i + \ell^y_j \geq lb^y_i + lb^y_j$ to this to get
		\[
			c^y_i + \frac{1}{2}\ell^y_i - (c^y_j - \frac{1}{2}\ell^y_j) \geq lb^y_i + lb^y_j,
		\]
		the desired inequality.
		\item \underline{$z^y_{i,j}=1, z^y_{j,i}=1$} Not feasible by \eqref{eqn:refined-proof-helper}.
	\end{itemize}
\endproof

\section{Alternative two-bit formulations}
\subsection{BLDP1 formulation} \label{app:bldp1}
When using the alternative codes $BB^4$, we may construct a formulation that is quite similar to the one presented in \eqref{eqn:gray-binary-formulation}, but with slightly different big-$M$ terms:

\begin{subequations}
\begin{alignat}{3}
	\frac{1}{2}\ell^s_k \leq c^s_k &\leq L^s - \frac{1}{2}\ell^s_k &\quad& \forall s \in \{x,y\}, k \in \{i,j\} \label{eqn:bb-first} \\
    \ell^s_k &\geq lb^s_k &\quad& \forall s \in \{x,y\}, k \in \{i,j\} \\
	c^y_i - c^y_j + \frac{1}{2}\left(\ell^y_i+\ell^y_j\right) &\leq L^y(y_1+y_2) \\
	c^x_i - c^x_j + \frac{1}{2}\left(\ell^x_i+\ell^x_j\right) &\leq L^x(2-y_1-y_2) \\
	c^y_j - c^y_i + \frac{1}{2}\left(\ell^y_i+\ell^y_j\right) &\leq L^y(1-y_1+y_2) \\
	c^x_j - c^x_i + \frac{1}{2}\left(\ell^x_i+\ell^x_j\right) &\leq L^x(1+y_1-y_2) \\
	y &\in \{0,1\}^2. \label{eqn:bb-last}
\end{alignat}
\end{subequations}
This formulation, with the addition of the area constraints \eqref{eqn:area}, forms the basis for the BLDP1 formulation from \cite{Castillo:2005}.

\subsection{Sequence-Pair formulation} \label{app:sequence-pair}
The sequence-pair formulation FLP-SP from \cite{Meller:2007} may be constructed from \eqref{eqn:gray-binary-formulation} with the addition of global constraints on the 0/1 variables, based on observations made by \cite{Murata:1996}. In particular, consider an $N$ box instance of the FLP and the corresponding formulation derived from Proposition \ref{prop:pairwise-to-Nbox} where each pairwise formulation $F^{i,j}$ is given by the gray binary formulation \eqref{eqn:gray-binary-formulation}. Then the addition of the following constraints yields the FLP-SP formulation:

\begin{subequations}
\begin{alignat}{2}
    \hat{w}^{i,j}_1 + \hat{w}^{j,k}_1 + \hat{w}^{k,i}_1 \leq 2 \quad &\forall i,j,k \in \llbracket N \rrbracket: i \neq j, i \neq k, j \neq k \\
    \hat{w}^{i,j}_2 + \hat{w}^{j,k}_2 + \hat{w}^{k,i}_2 \leq 2 \quad &\forall i,j,k \in \llbracket N \rrbracket: i \neq j, i \neq k, j \neq k,
\end{alignat}
\end{subequations}
where notationally
\[
    \hat{w}^{p,q}_k \defeq \begin{cases}
        w^{p,q}_k & p < q \\
        1-w^{p,q}_k & \text{o.w.}
    \end{cases} \quad \forall k \in \{1,2\}, p,q \in \llbracket N \rrbracket : p \neq q.
\]

\section{Proof of Proposition\ref{prop:map-from-refined}} \label{app:map-from-refined}
\proof{}
    We prove the second, as the first follows in the same way. Consider a feasible layout $(\hat{c},\hat{\ell}) \in \scrL_{i,j}$ and take the corresponding feasible codes $W \defeq \left\{w \in GB^4 : (\hat{c},\hat{\ell},w) \in \En(Q^{FLP},D^4,GB^4) \right\}$.
    Choose some $w \in W$. Then there exists some code $z \in C^8$ such that $(\hat{c},\hat{\ell},z) \in \En(Q^{FLP},D^8,C^8)$ and $\scrA^{GB}(w) \leq z$. Therefore, we have that, since $d \geq 0$,
    \[
        a^T\hat{c} + b^T\hat{\ell} + d^T\scrA^{GB}(w) \leq a^T\hat{c} + b^T\hat{\ell} + d^Tz \leq f.
    \]
    Therefore, the given inequality holds for $\En(Q^{FLP},D^4,GB^4)$.
\endproof

\section{Proof of Proposition~\ref{prop:map-to-refined}}\label{app:map-to-refined}
\proof{}
    Consider a feasible layout $(\hat{c},\hat{\ell}) \in \scrL_{i,j}$ and take the corresponding feasible codes $U \defeq \left\{u \in U^4 : (\hat{c},\hat{\ell},u) \in \En(Q^{FLP},D^4,U^4) \right\}$ and $Z \defeq \left\{z \in C^8 : (\hat{c},\hat{\ell},z) \in \En(Q^{FLP},D^8,C^8) \right\}$.
    From the construction of the codes $C^8$ discussed in Section \ref{sec:refined-disjunction}, for each $z \in Z$, there exist some $u,u' \in U$ such that $(z^y_{i,j},z^y_{j,i}) = (u^y_{i,j},u^y_{j,i})$ and $(z^x_{i,j},z^x_{j,i}) = (u'^x_{i,j},u'^x_{j,i})$ (if $z \in U^4$, then $u = u'$). For example, if $z = (1,1,0,0)$, we see that $u = (1,0,0,0)$ and $u' = (0,1,0,0)$. Therefore, if $d^x_{i,j}=d^x_{j,i}=0$, we have that
    \[
        a^T\hat{c} + b^T\hat{\ell} + d^Tz = a^T\hat{c} + b^T\hat{\ell} + d^Tu \leq f;
    \]
    if $d^y_{i,j}=d^y_{j,i}=0$, then the same inequality holds with $u'$ in place of $u$. Therefore, any inequality valid for $\En(Q,D^4,U^4)$ is valid for $\En(Q,D^8,C^8)$.
\endproof

\section{Proof of Proposition \ref{prop:ub-cuts}} \label{app:ub-cuts}
\proof{}
	We prove by enumerating the possible values for the components of $z$ having support over the constraint, noting in particular that $z^s_{p,q} = z^s_{q,p} = 1$ is always infeasible. Recall also that
	\[
		z^s_{p,q} = 1 \Longrightarrow c^s_p + \frac{1}{2}\ell^s_p \leq c^s_q - \frac{1}{2}\ell^s_q. \label{eqn:ast} \tag{*}
	\]
	\subsection{Inequality \eqref{eqn:ub-c1}}
	\begin{itemize}
		\item \underline{$z^s_{q,p}=0$} Sum the constraints $c^s_p \geq \frac{1}{2}\ell^s_p$ and $ub^s_q \geq \ell^s_q$ to get the constraint $c^s_p + ub^s_q \geq \frac{1}{2}\ell^s_p + \ell^s_q$.
		\item \underline{$z^s_{q,p}=1$} Using \eqref{eqn:ast} and rearranging gives
			\[
				c^s_p - c^s_q \geq \frac{1}{2}\ell^s_p + \frac{1}{2}\ell^s_q;
			\]
			adding the constraint $c^s_q \geq \frac{1}{2}\ell^s_q$ gives the desired result
			\[
				c^s_p \geq \frac{1}{2}\ell^s_p + \ell^s_q.
			\]
	\end{itemize}
	\subsection{Inequality \eqref{eqn:ub-c2}}
	\begin{itemize}
		\item \underline{$z^r_{p,q} = z^r_{q,p} = 0$}
			Want $\ell^s_p + \ell^s_q \leq L^s$. Since $z^s_{p,q} + z^s_{q,p} + z^r_{p,q} + z^r_{q,p} \geq 1$, we must have that either $z^s_{p,q} = 1$ or $z^s_{q,p} = 1$; w.l.o.g. choose the second. Then rearranging from \eqref{eqn:ast} gives
			\begin{align*}
				\frac{1}{2}\left(\ell^s_p + \ell^s_q\right) &\leq c^s_q - c^s_p \\
				&\leq \left(L^s-\frac{1}{2}\ell^s_q\right) - \left( \frac{1}{2}\ell^s_p \right) \quad \text{from \eqref{eqn:sitb}} \\
				\Longrightarrow \ell^s_p + \ell^s_q &\leq L^s.
			\end{align*}

		\item \underline{$z^r_{p,q} + z^r_{q,p} = 1$} Want that $ub^s_p + ub^s_q \geq \ell^s_p + \ell^s_q$, which follows immediately from the upper bounds on $\ell$.
	\end{itemize}
\endproof

\section{Proof of Proposition \ref{prop:obj-cuts}} \label{app:obj-cuts}
\proof{}
	We prove by enumerating the possible values for the components of $z$ having support over the constraint, noting in particular that $z^s_{p,q} = z^s_{q,p} = 1$ is always infeasible.
	\subsection{Inequality \eqref{eqn:obj1}}
	\begin{itemize}
		\item \underline{$z^s_{p,q} = 0, z^s_{q,p} = 0$}
			Want to show that
			\[
				d^s_{i,j} \geq \frac{1}{2}(\ell^s_p + \ell^s_q) - L^s,
			\]
			but since $d^s_{i,j} \geq 0$ necessarily (sum \eqref{eqn:linearized-objective}) it suffices to show that $\ell^s_p + \ell^s_q \leq 2L^s$, which follows immediately from summing \eqref{eqn:sitb} constraints $c^s_p \leq L^s - \frac{1}{2}\ell^s_p$ with $\frac{1}{2}\ell^s_p \leq c^s_p$ to get that $\ell^s_p \leq L^s$. Applying this also for $p$ and summing the resulting inequalities gives the result.
		\item \underline{$z^s_{p,q} = 1, z^s_{q,p} = 0$}
			Want to show that
			\[
				d^s_{i,j} \geq \frac{1}{2}(\ell^s_p + \ell^s_q),
			\]
			We have from \eqref{eqn:ast} that $c^s_p + \frac{1}{2}\ell^s_p \leq c^s_q - \frac{1}{2}\ell^s_h$; adding the appropriate constraint in \eqref{eqn:linearized-objective} to this gives the result.
		\item \underline{$z^s_{p,q} = 0, z^s_{q,p} = 1$} Same argument as the previous case.
	\end{itemize}
	\subsection{Inequality \eqref{eqn:obj2}}
	\begin{itemize}
		\item \underline{$z^s_{p,q} = 0, z^s_{q,p} = 0$}
			Want to show that
			\[
				d^s_{i,j} \geq c^s_p - c^s_h + \ell^s_p - L^s,
			\]
			which follows immediately from \eqref{eqn:linearized-objective} and the fact that $L^s \geq \ell^s_p$ for any feasible solution.
		\item \underline{$z^s_{p,q} = 1, z^s_{q,p} = 0$}
			Want to show that
			\[
				d^s_{i,j} \geq c^s_p - c^s_q + \ell^s_p + lb^s_q.
			\]
			Rearranging \eqref{eqn:ast} gives
			\[
				c^s_p - c^s_q + \frac{1}{2}\ell^s_p + \frac{1}{2}\ell^s_q \leq 0.
			\]
			Now summing the relation from \eqref{eqn:ast} with one of \eqref{eqn:linearized-objective} gives
			\[
				\frac{1}{2}(\ell^s_p + \ell^s_q) \leq d^s_{i,j}.
			\]
			Summing these two derived inequalities with $lb^s_q \leq \ell^s_q$ gives the desired result.
		\item \underline{$z^s_{p,q} = 0, z^s_{q,p} = 1$}
			Want to show that
			\[
				d^s_{i,j} \geq c^s_p - c^s_q + \ell^s_p + lb^s_q - L^s.
			\]
			Take the sum of one of \eqref{eqn:linearized-objective} and the inequality from \eqref{eqn:ast} to derive
			\[
				d^s_{i,j} \geq \frac{1}{2}\ell^s_p + \frac{1}{2}\ell^s_q.
			\]
			Furthermore, using our big-$M$ value, we have
			\[
				L^s \geq c^s_p - c^s_q + \frac{1}{2}(\ell^s_p + \ell^s_q);
			\]
			summing the two derived inequalities along with the lower bounds on $\ell$ gives the result.
	\end{itemize}
	\subsection{Inequality \eqref{eqn:obj3}}
	\begin{itemize}
		\item \underline{$z^s_{p,q} = 0$}
			Want to show that $d^s_{i,j} \geq c^s_p - c^s_q$, which is immediate from \eqref{eqn:linearized-objective}.
		\item \underline{$z^s_{p,q} = 1$}
			Want to show that
			\[
				d^s_{i,j} \geq c^s_p - c^s_q + lb^s_p + lb^s_q.
			\]
			From an argument above,
			\[
				d^s_{i,j} \geq \frac{1}{2}\ell^s_p + \frac{1}{2}\ell^s_q
			\]
			for this particular setting, and so we are done by summing this with
			\[
				c^s_p + \frac{1}{2}\ell^s_p \leq c^s_q - \frac{1}{2}\ell^s_q
			\]
			implied by \eqref{eqn:ast} and using the lower bounds on $\ell$.
	\end{itemize}
	\subsection{Inequality \eqref{eqn:obj4}}
	\begin{itemize}
		\item \underline{$z^s_{p,q} = 0, z^s_{q,p} = 0$}
			Want to show that
			\[
				2d^s_{i,j} \geq \ell^s_p - L^s,
			\]
			which just follows from the fact that $d^s_{i,j} \geq 0$ and $L^s \geq \ell^s_p$.
		\item \underline{$z^s_{p,q} = 1, z^s_{q,p} = 0$}
			Want to show that
			\[
				d^s_{i,j} \geq \frac{1}{2}\ell^s_p + \frac{1}{2}lb^s_q,
			\]
			which follows immediately from the inequality (valid for this particular setting for $z$)
			\[
				d^s_{i,j} \geq \frac{1}{2}\ell^s_p + \frac{1}{2}\ell^s_q.
			\]
			derived previously.
		\item \underline{$z^s_{p,q} = 0, z^s_{q,p} = 1$} Same argument as the previous case.
	\end{itemize}
\endproof

\section{Proof of Proposition~\ref{prop:multibox-cuts}} \label{app:multibox-cuts}
\proof{}
    First, we note that inequalities (\ref{eqn:multi1}-\ref{eqn:multi7}) are variations on (\ref{eqn:obj1}-\ref{eqn:obj4},\ref{eqn:ru-tight-sitb},\ref{eqn:ru-tight-sitb},\ref{eqn:ru-nonoverlap}), respectively, with an additional $\gamma_P\scrM^s_P(z)$ term appearing.

	We can use \eqref{eqn:z-variable-implication} to see that $z^s_{t^{\xi-1},t^\xi} = 1$ implies $c^s_{t^{\xi-1}} + \frac{1}{2}\ell^s_{t^{\xi-1}} \leq c^s_{t^{\xi}} - \frac{1}{2}\ell^s_{t^{\xi}}$ for all $\xi \in \llbracket m+1 \rrbracket$ to take a telescoping sum and derive $c^s_{j} - \frac{1}{2}\ell^s_{j} \geq c^s_{i} + \frac{1}{2}\ell^s_{i} + \sum_{\xi=1}^m \ell^s_{t^\xi} \geq c^s_{i} + \frac{1}{2}\ell^s_{i} + \gamma_P$
	in the case where $z^s_{t^{\xi-1},t^\xi} = 1$ for all $\xi \in \llbracket m+1 \rrbracket$. Combining this with property \eqref{eqn:multibox-map-property}, we derive that
	\[
	    c^s_{j} - \frac{1}{2}\ell^s_{j} \geq c^s_{i} + \frac{1}{2}\ell^s_{i} + \gamma_P\scrM(z) \tag{**}
	\]
	is valid for any feasible solution for $F^{RU}$. This can be used to directly derive (\ref{eqn:multi5}-\ref{eqn:multi7}). For example, \eqref{eqn:multi5} can be derived by summing the valid inequalities
	\[
	    c^s_{i} + \frac{1}{2}\ell^s_{i} + \gamma_P\scrM_P(z) \leq c^s_{j} - \frac{1}{2}\ell^s_{j}, \quad \quad
	    \frac{1}{2}\ell^s_i \leq c^s_i, \quad \quad
	    lb^s_iz^s_{i,j} \leq \ell^s_i,
    \]
    after noting that $z^s_{i,j} \in \{0,1\}$ implies that $lb^s_{i,j}z^s_{i,j} \leq lb^s_{i,j}$.

	For (\ref{eqn:multi1}-\ref{eqn:multi4}), we first observe that, due to property \eqref{eqn:multibox-map-property}, the case analysis in Appendix~\ref{app:obj-cuts} will only differ in the case where $z^s_{t^{\xi-1},t^\xi} = 1$ for all $\xi \in \llbracket m+1 \rrbracket$. However, under the assumption that $lb^s_k > 0$ for all $k \in \llbracket N \rrbracket$, this also implies that $z^s_{i,j} = 1$. Therefore, we may use the tightened inequality (**) in lieu of (*) in the case analyses, yielding the result.
\endproof

\section{Valid inequalities from the literature} \label{app:literature-cuts}
We now present the B2 and V2 inequalities from \cite{Meller:1999} in the notation used in the present work. The B2 inequalities for the unary formulation are of the form
\begin{equation}
    d^s_{i,j} \geq \frac{1}{2}(lb^s_i+lb^s_j)(u^s_{i,j}+u^s_{j,i}) \quad \forall s \in \{x,y\}.
\end{equation}
The V2 inequalities for the unary formulation are
\begin{equation}
    d^s_{i,j} \geq \frac{1}{2}(\ell^s_i+\ell^s_j) - \frac{1}{2}\min\{ub^s_i+ub^s_j,L^s\}(1-u^s_{i,j}-u^s_{j,i}) \quad \forall s \in \{x,y\},
\end{equation}
which is equivalent \eqref{eqn:obj1} with a potentially tightened coefficient $\frac{1}{2}\min\{ub^s_i+ub^s_j,L^s\}$.
Using Proposition \ref{prop:map-to-refined} and Proposition \ref{prop:map-from-refined}, these inequalities may be applied to all the formulations discussed in this work.

\section{Symmetry-breaking} \label{app:symmetry-breaking}
The symmetry-breaking described in \cite{Sherali:2003} works by restricting the possible relative layout between a single pair of components in a modification of the so-called \emph{position $p-q$ method} from \cite{Meller:1999}. The scheme chooses a single pair $(p,q) \in \scrP$; in this work, we follow \cite{Sherali:2003} and choose $(p,q) \in \arg\max_{(i,j) \in \scrP} p_{i,j}$. We then may add the following constraints to the refined unary formulation:
\begin{subequations}
\begin{alignat}{2}
    c^s_p &\leq c^s_q &\quad \forall s \in \{x,y\} \\
    z^s_{q,p} &\leq 0 &\quad \forall s \in \{x,y\} \\
    (c^x_q - c^x_p) + (c^y_q - c^y_p) &\geq \frac{1}{2}\min\{lb^x_p+lb^x_q,lb^y_p+lb^y_q\}&.
\end{alignat}
\end{subequations}
Using Proposition \ref{prop:map-from-refined}, these inequalities may be applied to all the formulations discussed in this work.

\section{Relative root and relaxation gap (Tables)} \label{app:relaxation-gap-table}
In this section, we present results for a single instance \texttt{base\_instance-0.0(5)} for each instance family \texttt{base\_instance}, as the other instances do not provide additional insight. for brevity, we will omit the perturbation and aspect ratio factors from the names.
\begin{table}
    {\footnotesize
    \centering
    \begin{tabular}{r|ccc}
                              &   \#1 & \#1 w/ SB & \#2 \\ \hline
        \texttt{hp}           & 100\% &            89.0\% &    51.5\% \\
        \texttt{apte}         & 100\% &            87.5\% &    58.4\% \\
        \texttt{xerox}        & 100\% &            84.6\% &    56.2\% \\
        \texttt{Camp91}       & 100\% &            77.1\% &    43.8\% \\
        \texttt{Bozer97-1}    & 100\% &            88.7\% &    61.6\% \\
        \texttt{Bozer97-2}    & 100\% &            93.6\% &    54.6\% \\
        \texttt{Bazaraa75-1}  & 100\% &            85.8\% &    63.1\% \\
        \texttt{Bazaraa75-2}  & 100\% &            87.8\% &    68.7\% \\
        \texttt{Bozer91}      & 100\% &            94.4\% &    43.7\% \\
        \texttt{Armour62-1}   & 100\% &            96.5\% &    65.8\% \\
        \texttt{Armour62-2}   & 100\% &            96.8\% &    64.3\%
    \end{tabular}
    \caption{Relative gap of the relaxation lower bound, with respect to the best known feasible solution. Group \#1 includes \texttt{U}, \texttt{BLDP1}, \texttt{BLDP1+}, \texttt{SP}, \texttt{SP+}, \texttt{SP+VI}, \texttt{SP+VI3}, and \texttt{RU}. Group \#2 includes \texttt{U+}, \texttt{RU+VI}, and \texttt{RU+VI3}. Symmetry breaking from \cite{Sherali:2003} is added to Group \#1 for comparison (\# w/ SB); it does not affect the values for Group \#2.}
    }
    \label{table:relaxation}
\end{table}

\begin{table}
    \centering
    {\footnotesize
    \begin{tabular}{r|cccccc}
                & \texttt{U}/\texttt{RU} & \texttt{U+}/\texttt{U+VI}/\texttt{RU+VI} & \texttt{BLDP1} & \texttt{SP}/\texttt{SP+} & \texttt{SP+VI} & \texttt{SP+VI3} \\ \hline
        \texttt{hp}   & 69.9\% &                     51.5\% & 89.0\% & 88.4\% & 83.0\% & 83.7\% \\
        \texttt{apte} & 72.9\% &                     58.4\% & 85.7\% & 86.4\% & 83.6\% & 85.2\% \\
        \texttt{xerox} & 70.3\% &                     56.2\% & 84.6\% & 83.6\% & 79.1\% & 78.4\% \\
        \texttt{Camp91} & 60.4\% &                    40.2\% & 77.1\% & 73.3\% & 48.9\% & 52.6\% \\
        \texttt{Bozer97-1} & 74.8\% &                 61.0\% & 85.5\% & 85.7\% & 79.7\% & 77.8\% \\
        \texttt{Bozer97-2} & 73.2\% &                 50.9\% & 93.6\% & 93.2\% & 73.9\% & 76.3\% \\
        \texttt{Bazaraa75-1} & 72.6\% &               61.3\% & 85.2\% & 84.9\% & 77.9\% & 78.3\% \\
        \texttt{Bazaraa75-2} & 78.1\% &               68.7\% & 85.0\% & 87.5\% & 82.8\% & 83.1\% \\
        \texttt{Bozer91} & 69.1\% &                   42.2\% & 94.4\% & 94.4\% & 61.7\% & 59.9\% \\
        \texttt{Armour62-1} & 81.2\% &                65.8\% & 96.5\% & 95.3\% & 93.3\% & 91.3\% \\
        \texttt{Armour62-2} & 80.5\% &                64.3\% & 96.8\% & 96.6\% & 92.0\% & 94.7\%
    \end{tabular}
    \caption{Relative gap of the root node lower bound, with respect to the best known feasible solution. }
    }
    \label{table:root}
\end{table}

\end{document}